\documentclass[fontsize=12pt,a4paper,headings=normal,
twoside=false,leqno,parskip=half-,abstract=true]{scrartcl}
\usepackage[english]{babel}
\usepackage[utf8]{inputenc}
\usepackage{hyperref}
\hypersetup{
 pdftitle={Unbounded Sturm Attractor},
 pdfauthor={Phillipo Lappicy, Juliana Pimentel},
 colorlinks=true,
 linkcolor=blue,
 citecolor=blue,
 filecolor=blue,
 urlcolor=blue}

\usepackage{graphicx}
\usepackage[format=plain,labelfont=bf,font=small]{caption}
\usepackage{subfigure}
\usepackage{xcolor}
\usepackage[arrow, matrix, curve]{xy}
\usepackage{float}

\usepackage{caption}
\usepackage{tabulary}
\usepackage{array}
\newcolumntype{N}[1]{>{\centering\arraybackslash}m{#1}}

\usepackage{amsmath,amsthm}
\usepackage{amssymb} 

\makeatletter
\newcommand{\tpitchfork}{%
  \vbox{
    \baselineskip\z@skip
    \lineskip-.52ex
    \lineskiplimit\maxdimen
    \m@th
    \ialign{##\crcr\hidewidth\smash{$-$}\hidewidth\crcr$\pitchfork$\crcr}
  }%
}
\makeatother
\usepackage{latexsym}

\usepackage[notref,notcite,color,final 
]{showkeys}

\definecolor{refkey}{rgb}{1,0,0}
\definecolor{labelkey}{rgb}{1,0,0}

\usepackage{tikz}
\usepackage{tikz-cd}

\usepackage[textwidth=2cm,textsize=small,backgroundcolor=none]{todonotes}

  \mathchardef\ordinarycolon\mathcode`\:
  \mathcode`\:=\string"8000
  \begingroup \catcode`\:=\active
    \gdef:{\mathrel{\mathop\ordinarycolon}}
  \endgroup

\newtheorem{thm}{Theorem}[section]
\newtheorem{lem}[thm]{Lemma}
\newtheorem{prop}[thm]{Proposition}

\newtheorem{cor}[thm]{Corollary}

\usepackage{comment}
\usepackage{caption}

\begin{document}

\title{{\LARGE{Unbounded Sturm attractors for \\ quasilinear parabolic equations}}}

\author{
 \\
{~}\\
Phillipo Lappicy* and Juliana Fernandes**\\
\vspace{2cm}}

\date{ }
\maketitle
\thispagestyle{empty}

\vfill

$\ast$\\
Instituto de Ciências Matemáticas e de Computação, Universidade de S\~ao Paulo\\
Av. trabalhador são-carlense, São Carlos, Brazil\\
$\ast$\\
Depto. Análisis Matemático y Mat. Aplicada, Universidad Complutense de Madrid\\
Pl. de las Ciencas 3, Madrid, Spain\\
\\
$\ast \ast$\\
Instituto de Matemática, Universidade Federal do Rio de Janeiro\\
Av. Athos da Silveira Ramos, Rio de Janeiro, Brazil


\newpage
\pagestyle{plain}
\pagenumbering{arabic}
\setcounter{page}{1}

\begin{abstract}
\noindent
We analyze the asymptotic dynamics of quasilinear parabolic equations when solutions may grow up (i.e., blow up in infinite time). For such models, there is a global attractor which is unbounded and the semiflow induces a nonlinear dynamics at infinity by means of a Poincaré projection. 
In case the dynamics at infinity is given by a semilinear equation, then it is gradient, consisting of the so-called equilibria at infinity and their corresponding heteroclinics. Moreover, the diffusion and reaction compete for the dimensionality of the induced dynamics at infinity.
If the equilibria are hyperbolic, we explicitly prove the occurrence of heteroclinics between bounded equilibria and/or equilibria at infinity.
These unbounded global attractors describe the space of admissible initial data at event horizons of certain black holes.

\ 

\end{abstract}

\section{Main results}\label{sec:intro}

\numberwithin{equation}{section}
\numberwithin{figure}{section}
\numberwithin{table}{section}

Consider the scalar quasilinear parabolic differential equation
\begin{equation}\label{PDE}
    u_t = \mathcal{L}(u):= a(x,u,u_x) u_{xx}+bu+f(x,u,u_x)
\end{equation}
with initial data $u(0,x)=u_0(x)$ such that $a,f:[0,\pi]\times \mathbb{R}^2\rightarrow \mathbb{R}$ are bounded $C^2$ functions satisfying the strict parabolicity condition $a(x,u,u_x)\geq \epsilon>0$, and $x \in [0,\pi]$ with Neumann boundary conditions. 

The equation \eqref{PDE} defines a semiflow, denoted by $(t,u_0)\mapsto u(t)$, in a Banach space $X^\alpha:=C^{2\alpha+\beta}([0,\pi])\cap \{\text{Neumann b.c.}\}$, which is a space of Hölder continuous functions to be defined in Section \ref{sec:func}. We suppose $2\alpha+\beta>1$, so that solutions are in $C^1([0,\pi])$. 

The class of equations \eqref{PDE} that only exhibit bounded solutions is called \emph{dissipative}. In this setting, if $b=0$ and $f$ satisfy certain growth conditions, there exists a global attractor, $\mathcal{A}\subseteq X^\alpha$, given by the maximal compact invariant set that attracts all bounded sets, see \cite{BabinVishik92}. 
Moreover, the dynamics in $\mathcal{A}$ is gradient, due to a Lyapunov function constructed by Zelenyak and Matano \cite{Zelenyak68,Matano88} given by 
\begin{equation}\label{BOOK}
    E:= \int_{0}^\pi L(x,u,u_x) dx \quad \text{ such that } \quad   \frac{dE}{dt}= -\int_{0}^\pi \frac{|u_t|^2}{a(x,u,u_x)} dx\leq 0.
\end{equation}
Therefore the attractor is decomposed as $\mathcal{A}=\mathcal{E}\cup \mathcal{H}$, where $\mathcal{E}$ denotes the equilibria points (time independent solutions, i.e., $u_t=0$) and $\mathcal{H}$ stands for the set of heteroclinic orbits, i.e., a solution $u(t)\in \mathcal{H}$ satisfies
\begin{equation}\label{defHETS}
        e_j\xleftarrow{t\to -\infty}u(t)\xrightarrow{t\to \infty} e_k, \qquad \text{ where } e_j,e_k\in\mathcal{E}, j\neq k.
\end{equation}
The task of explicitly finding equilibria and which heteroclinics occur is called the connection problem.
In particular, necessary and sufficient conditions are given in order to guarantee the occurrence of heteroclinics among two given hyperbolic equilibria as in \eqref{defHETS}.
Such construction of the global attractor was carried out in the semilinear context ($a\equiv 1$) with Hamiltonian reaction $f(u)$ by Brunovský and Fiedler \cite{FiedlerBrunovsky89}, whereas the more general reaction $f(x,u,u_x)$ was treated by Fiedler and Rocha \cite{FiedlerRocha96}, and in case of periodic boundary conditions by Fiedler, Rocha and Wolfrum \cite{FRW}. Quasilinear and fully nonlinear equations were pursued by one of the authors \cite{LappicyQuasi,LappicyFully}. 
These attractors are known as \emph{Sturm attractors} because the connection problem can be solved by means of nodal properties discovered by Sturm \cite{Sturm}, and generalized by Matano \cite{Matano82} and Angenent \cite{Angenent88}.

On the other hand, the class of equations \eqref{PDE} that admit unbounded solutions are called \emph{non-dissipative}. 
Solutions $u(t)$ that become unbounded in $X^\alpha$ as $t\to T$ can be divided into two categories. First, \emph{finite time blow-up solutions},  when $T<\infty$, 
are part of an established and active research topic, see \cite{FiMa02} and references therein. 
Second, \emph{infinite time blow-up solutions} (also called \emph{grow-up solutions}), when $T=\infty$, 
have attracted a more recent attention.
For the existence of grow-up in 
equations with a localized reaction, see \cite{ArrietaPardo,FerreiradePablo20}, and in case of a fractional diffusion, see \cite{MussoSire19}.
In the non-dissipative setting, there does not exist a global attractor in the usual sense, which is the maximal compact invariant set, see \cite{Henry81}. Yet, there is an \emph{unbounded global attractor} $\mathcal{A}\subseteq X^\alpha$ which is defined as the minimal invariant non-empty set in $X^\alpha$ attracting all bounded sets, firstly introduced by Chepyzhov and Goritskii \cite{CG92}. See also \cite{BortolanJu22}.

We desire to describe the dynamical behavior of unbounded solutions of \eqref{PDE}. 
In Lemma \ref{lemnondiss}, we will show that $b>0$ is a sufficient condition for the existence of grow-up solutions.
Moreover, we will decompose the unbounded attractor $\mathcal{A}$ into smaller invariant sets, describe them and show how they are related by means of heteroclinics, as in the dissipative case.
Such attractors are known as \emph{unbounded Sturm attractors}. 

Despite non-dissipativity, the parabolic equation \eqref{PDE} still possess a Lyapunov function according to \eqref{BOOK}, as long as solutions exist.
Therefore, in case of hyperbolic equilibria, the following dichotomy holds: either a solution remains bounded and converges to a bounded equilibrium as $t\to\infty$, or it is a grow-up solution. See \cite[Section 4.3]{Henry81} and \cite{QuittnerSouplet19}.
The grow-up solutions were interpreted as heteroclinic orbits to infinity by Hell in \cite{Hell09}. In order to describe the dynamics of unbounded solutions and account for \emph{how} solutions grow-up, an infinite dimensional sphere $\mathcal{S}^\infty$ was added at infinity with an induced semiflow, 
by means of a Poincaré projection. 
Previous investigations to understand such structure at infinity for semilinear equations were done by Hell \cite{Hell09}, Ben-Gal \cite{BenGal10} for $f(u)$, Pimentel and Rocha \cite{RochaPimentel16} for $f(x,u,u_x)$, and Pimentel \cite{P} for periodic boundary condition. 
We now describe this process in detail.

The Poincaré projection maps the phase-space $X^{\alpha}$ of \eqref{PDE} to a subset of the unit sphere in $L^2\times \mathbb{R}$. 
Indeed, identify $X^\alpha$ with $X^{\alpha}\times\{1\}\subseteq L^2\times\{1\}$. The set $L^2\times\{1\}$ is the tangent space (at the north pole) of the northern hemisphere, $\mathbb{S}_+\subseteq L^2\times \mathbb{R}$, called the \emph{Poincaré hemisphere} and given by
\begin{equation}\label{Phemis}
    \mathbb{S}_+:=\{(\chi,z)\in L^2\times [0,1]: ||\chi||^{2}_{L^2}+z^{2}=1\}.
\end{equation}
Then for each point in phase-space, $u\in X^{\alpha}\subseteq L^2$, consider the line that passes through $(u,1)\in L^2\times [0,1]$ and the origin $(0,0)\in L^2\times [0,1]$. This line intersects the upper hemisphere $\mathbb{S}_+$ at a point, which defines the projection $\mathcal{P}:X^\alpha\times\{1\}\to \mathbb{S}_+$, called the \emph{Poincaré projection}\footnote{The Poincaré \emph{projection} is often called \emph{compactification}, since it compactifies finite dimensional spaces. However, it does not compactify the infinite dimensional phase-space $X^\alpha$. 
}. See Figure~\ref{Fig:P}.
%
\begin{figure}[H]\centering
    \begin{tikzpicture}[scale=2]
    \draw[->] (-1.5,0) -- (1.5,0) node[right] {$\chi\in L^2 $};
    \draw[->] (0,-0.1) -- (0,1.5) node[above] {$z\in \mathbb{R}$};

    \draw [thick, domain=0:3.14,variable=\t,smooth] plot ({cos(\t r)},{sin(\t r)}); 
    
    \draw[gray!75,-] (1.5,1) -- (-1.5,1) node[left] {$X^\alpha\hookrightarrow L^2\times \{ 1\}$};

    \filldraw [black] (1,1) circle (1pt) node[anchor=south]{$(u,1)$};
    \filldraw [black] (0,0) circle (1pt) node[anchor=north]{$(0,0)$};
    \filldraw [black] (0.7,0.7) circle (1pt) node[anchor=west]{$\mathcal{P}(u,1)$};
    \filldraw [black] (-0.707,0.707) circle (0.1pt) node[anchor=east]{$\mathbb{S}_+$};
    \filldraw [black] (1,0) circle (1pt) node[anchor=north]{$\mathcal{S}^\infty$};
    \filldraw [black] (-1,0) circle (1pt);

    \draw[-] (0,0) -- (1,1);
    \end{tikzpicture}
\caption{Poincaré projection $\mathcal{P}$ from phase-space $X^\alpha\hookrightarrow L^2\times \{ 1\}$ into the hemisphere $\mathbb{S}_+$. As solutions grow-up, $||u(t)||_{L^2}\to \infty$, the projection $\mathcal{P}(u,1)$ converges to the equator $\mathcal{S}^\infty:=\mathbb{S}_+|_{z=0}$.}\label{Fig:P}
\end{figure}
The coordinates of the projection $\mathcal{P}(u,1)$ are denoted by $(\chi,z)$ and can be computed from the colinearity of the points $(0,0)$, $(u,1)$, and the intersection at $(\chi,z)$ with the hemisphere $\mathbb{S}_+$, yielding 
\begin{equation}\label{P}
    (\chi,z):={\mathcal{P}}(u,1)=\left(\frac{u}{ \sqrt{1+||u||_{L^2}^2}},\frac{1}{ \sqrt{1+||u||_{L^2}^2}}\right).
\end{equation}
%
Therefore, Hell's perspective that grow-up solutions in $X^{\alpha}$ are heteroclinics to infinity can be interpreted (in terms of the Poincaré projection $\mathcal{P}$) as solutions in the hemisphere $\mathbb{S}_+$ that converge to the equator, which is characterized by $z=0$ and $||\chi||^2_{L^2}=1$. Indeed, note that $z=1$ if, and only if, $u \equiv 0$. Hence, the origin of $X^{\alpha}$ is mapped to the north pole of $\mathbb{S}_+$. Moreover, $z$ decreases to $0$ if, and only if, $||u||_{L^2}$ increases to $\infty$. 
Therefore the relevant asymptotic unbounded dynamics of the projected semiflow is contained in the unit sphere of $L^2$, consisting of bounded trajectories with coordinates $(\chi,0)$ such that $||\chi||^2_{L^2}=1$. For this reason, the equator of $\mathbb{S}_+$ is called the \textit{sphere at infinity}, and it is denoted by 
\begin{equation}
    \mathcal{S}^\infty:=\mathbb{S}_+|_{z=0}.
\end{equation}
The projection $\mathcal{P}$ induces a semiflow on $\mathbb{S}_+|_{z>0}$, which is obtained by a homothety (with scale factor $z$) of the original vector field $\mathcal{L}(u)$ in \eqref{PDE}. 
%
Indeed, differentiating \eqref{P} with respect to time, the new variables $(\chi,z)$ satisfy
\begin{subequations}\label{flowSPHERE}
\begin{align}
    \chi_t&
    = \mathcal{L}^z(\chi)-\chi\langle \mathcal{L}^z(\chi),\chi\rangle  \label{flowSPHERE1}\\
    z_t& 
    =-\langle \mathcal{L}^z(\chi),\chi\rangle \cdot z\label{flowSPHERE2}
\end{align}
where the projected vector field depends on a homothety of the original vector field \eqref{PDE} with scale factor $z:=(1+||u||^2)^{-1/2}$ defined by
\begin{equation}\label{homotvf}
    \mathcal{L}^z(\chi):=z \mathcal{L}(z^{-1}\chi)=a^z(x,\chi,\chi_x)\chi_{xx}+b\chi+f^z(x,\chi,\chi_x),
\end{equation}
\end{subequations}
where $a^z(x,\chi,\chi_x):=a(x,z^{-1}\chi, z^{-1}\chi_x)$ and $f^z(x,\chi,\chi_x):=z f(x, z^{-1}\chi, z^{-1}\chi_x)$ are homotheties (with the same scale factor $z>0$) of $a$ and $f$, respectively. 

The projection \eqref{P} thereby induces a semiflow within $\mathbb{S}_+|_{z>0}$ described by equations \eqref{flowSPHERE}.
Moreover, the induced (nonlinear and nonlocal) semiflow at $\mathcal{S}^\infty$, i.e. for $z=0$, is given by the limit as $z\to 0$, which may produce, in a number of settings, a degenerate or singular semiflow at the sphere at infinity $\mathcal{S}^\infty$. See \cite{BenGal10,Hell09} for further details and examples.
%
For this reason, we restrict to the case that the projected quasilinear diffusion coefficient converges uniformly, so that there is a well defined limiting semiflow at the invariant subspace $\mathcal{S}^\infty$. Mathematically, we suppose the following asymptotic condition for the real-valued quasilinear diffusion coefficient:
\begin{equation}\label{alim}
    \lim_{|(u,p)|_{\mathbb{R}^2}\to \infty} a(x,u,p)=a^\infty,
\end{equation}
for all $x\in [0,\pi]$ and some $a^\infty\in\mathbb{R}_+$. Note that \eqref{alim} denotes a limit in $\mathbb{R}$ for any fixed $x\in [0,\pi]$ and any direction $(u,p)\in \mathbb{R}^2$ that goes to infinity.
Note the convergence in $\mathbb{R}$ implies the pointwise convergence in $L^2$ of its associated Nemitskii operator to $a^\infty$. In particular, the homothety $a^z(x,\chi, \chi_x)$ converges pointwise (in $L^2$) to $a^\infty$ as $z\to 0$.

Therefore, the Poincaré projection \eqref{P} transforms unbounded grow-up solutions $u(t)$ of \eqref{PDE} in $X^\alpha$ into bounded solutions $(\chi(t),z(t))$ of \eqref{flowSPHERE} in the hemisphere $\mathbb{S}_+$ that converge to the invariant equator: the sphere at infinity $\mathcal{S}^\infty$. See \cite{Hell09,BenGal10}.
We are interested in describing the dynamics of the projected semiflow that corresponds to the asymptotic unbounded dynamics of $X^\alpha$. In particular, the unbounded global attractor $\mathcal{A}\subseteq X^\alpha$ is projected into the sphere, $\mathcal{P}(\mathcal{A})\subseteq \mathbb{S}_+$, which is encapsulated by a subset in the sphere at infinity $\mathcal{S}^\infty$. We refer to the closure (in $\mathbb{S}_+$) of $\mathcal{P}(\mathcal{A})$ as the \emph{extended unbounded attractor} and denote it by 
\begin{equation}\label{extendedA}
    \overline{\mathcal{P}(\mathcal{A})}:=\mathrm{cl}_{\mathbb{S}_+}(\mathcal{P}(\mathcal{A})),    
\end{equation}
which includes the limiting unbounded dynamics at infinity.

Below we present the first main result. We prove that the dynamics induced in the closed hemisphere $\mathbb{S}_+$ is gradient.
Besides the conclusion regarding the bounded dynamics, which consists of (bounded) equilibria and (bounded) heteroclinics between them, there are two conclusions to be drawn regarding the unbounded dynamics. 

First, the $\omega$-limit of a projected grow-up solution consists of a single equilibrium of the induced semiflow at $\mathcal{S}^\infty$, which is called \emph{equilibrium at infinity}\footnote{Alternatively, one can consider the extension of an unbounded solution by infinity, which is called a \emph{metasolution}, see \cite{LopezGomez03,LopezGomez16} and references therein. However, metasolutions are objects which do not belong to phase-space, in contrast to equilibria at infinity, which lie in $\mathcal{S}^\infty$. 
The profiles of equilibria at infinity are compactified profiles of metasolutions through \eqref{P}.
For the approximation of metasolutions by solutions, see \cite{LopezGomez01}, which can be compared to the approximation scheme in \cite{BCP} using our current Poincaré projection approach.
} and denoted by $\pm\Phi_{j}\in \mathcal{S}^\infty$.
This occurs because unbounded solutions of \eqref{PDE} grow-up more rapidly in the direction of certain eigenfunctions of the operator $- a^\infty \partial^2_x$, given by $\{\varphi_j\}_{j\in\mathbb{N}_0}$, and thereby the equilibria of the induced flow at infinity are explicitly given by
\begin{equation}\label{CptEq}
    \pm\Phi_{j}:=(\pm\varphi_{j},0) \in \mathcal{S}^\infty.
\end{equation}
Second, the remaining relevant dynamics at $\mathcal{S}^\infty$ are heteroclinics between equilibria at infinity. 
Thus, no complicated dynamics arises in the extended attractor $\overline{\mathcal{P}(\mathcal{A})}$. 
%
%
%
\begin{thm}\emph{\textbf{Decomposition of the extended unbounded Sturm Attractor.} } \label{attractorthm1}
Let $a,f\in C^2$ be bounded such that $a$ is strictly parabolic with limiting behavior \eqref{alim}. Fix $b>0$ such that $\sqrt{b/a^\infty}\not\in \mathbb{N}$. Suppose bounded equilibria are hyperbolic. Then, the extended unbounded attractor $\overline{\mathcal{P}(\mathcal{A})}\subseteq \mathbb{S}_+$ of \eqref{PDE} is decomposed as 
    \begin{equation}
            \overline{\mathcal{P}(\mathcal{A})}=\overline{\mathcal{E}}\cup \overline{\mathcal{H}}, 
    \end{equation}
where the set of extended equilibria $\overline{\mathcal{E}}$ consists of:
\begin{itemize}
    \item[(i)] projected bounded equilibria, $\mathcal{E}^b:=\mathcal{P}(\mathcal{E})$, where $\mathcal{E}:=\{e_j: j=1,\dots, N\}$ is the set of bounded equilibria of \eqref{PDE} for some $N\in \mathbb{N}$,  
    \item[(ii)] equilibria at infinity, $\mathcal{E^\infty}=\{\pm\Phi_j:j=0,\dots, N^\infty\}$, where $N^\infty:= \lfloor\sqrt{b/a^\infty}\rfloor$,
\end{itemize}
whereas the set of extended heteroclinics $\overline{\mathcal{H}}$ consists of:
\begin{itemize}
    \item[(i)] projected bounded heteroclinic orbits, $\mathcal{H}^b:=\mathcal{P}(\mathcal{H})$, where $\mathcal{H}$ is the set of bounded heteroclinics of \eqref{PDE} ,
    \item[(ii)] projected grow-up solutions, $\mathcal{H}^{up}:=\mathcal{P}(\{u(t)\in X^\alpha : ||u||_\alpha\to \infty\})$, which can be seen as heteroclinics from bounded equilibria to equilibria at infinity,
    \item[(iii)] heteroclinics at infinity, $\mathcal{H^\infty}\subseteq \mathcal{S}^\infty$, between equilibria at infinity.
\end{itemize} 
\end{thm}
\begin{figure}[H]\centering
\vspace{-0.5cm}
\begin{subfigure}\centering
    \begin{tikzpicture}[scale=0.75]
    \draw[very thick, gray,-] (-2,0) -- (-1,0);
    \draw[very thick, gray,->] (-1.5,0) -- (-1.51,0) node[above] {\footnotesize{$\mathcal{H}^{up}$}};

    \draw[gray,-] (-1,0) -- (0,0);
    \draw[gray,<-] (-0.5,0) -- (-0.51,0) node[above] {\footnotesize{$\mathcal{H}^{b}$}};
        
    \filldraw [very thick] (0,-2) circle (2pt) node[below] {\footnotesize{$\mathcal{E}^\infty$}};
    \filldraw [very thick] (0,2) circle (2pt) node[above] {\footnotesize{$\mathcal{E}^\infty$}};    
    \filldraw [very thick] (-2,0) circle (2pt) node[left] {\footnotesize{$\mathcal{E}^\infty$}};
    \filldraw [very thick] (2,0) circle (2pt) node[right] {\footnotesize{$\mathcal{E}^\infty$}};
     
    \draw [very thick] (0,0) circle (57pt);

    \draw[very thick,<-] (1.2,1.6) arc (44:45:0.4cm and 0.4cm)  node[right] {\footnotesize{$\mathcal{H}^\infty$}};    
    \draw[rotate=180,thick,<-] (1.2,1.6) arc (44:45:0.4cm and 0.4cm);    

    \draw[rotate=70,very thick,->] (1.2,1.6) arc (44:45:0.4cm and 0.4cm);    
    \draw[rotate=70+180,very thick,->] (1.2,1.6) arc (44:45:0.4cm and 0.4cm);    

    \filldraw [gray,very thick] (-1,0) circle (1pt) node[anchor=north]{\footnotesize{$\mathcal{E}^b$}};
    \filldraw [gray,very thick] (0,0) circle (1pt)  node[anchor=north]{\footnotesize{$\mathcal{E}^b$}};
    \end{tikzpicture}    
\end{subfigure}
\vspace{-0.25cm}
\captionof{figure}{Schematic depiction of the extended unbounded attractor $\overline{\mathcal{P}(\mathcal{A})}\subseteq \mathbb{S}_+$. The bounded part consists of projected bounded equilibria (gray dots) and heteroclinics (gray arrow), $\mathcal{E}^b\cup \mathcal{H}^b$. The unbounded solutions consist of the projected grow-up solutions (thick gray arrow), $\mathcal{H}^{up}$, which are heteroclinics from projected bounded equilibria towards equilibra at infinity, $\mathcal{E}^\infty$. The extension of the projected attractor is the sphere at infinity, $\mathcal{S}^\infty$, which consists of equilibria at infinity (black dots) and their heteroclinics (black arrows), $\mathcal{E}^{\infty}\cup \mathcal{H}^{\infty}$.
}\label{FIG:scheme}
\end{figure}

Next, we give a detailed description of the structure at infinity. 
In particular, the dimension of such structure depends on the interplay between the asymptotic diffusion parameter $a^\infty\in\mathbb{R}_+$ and the linear reaction parameter $b>0$.
Moreover, we give necessary and sufficient conditions for the occurrence of bounded and/or unbounded heteroclinics.
This implies that only a finite-dimensional subset of the infinite-dimensional sphere at infinity is attainable by grow-up solutions. 

We present a few notions before the upcoming Theorem.
Denote by the \emph{zero number} $z(u_*)$ the number of sign changes of a function $u_*(x)$. 
Recall that an equilibrium $u_*$ is \emph{hyperbolic} if the linearization operator of the right hand side of \eqref{PDE} at $u_*$ has no eigenvalue being zero.
Also, the \emph{Morse index} $i(u_*)$ of a hyperbolic equilibrium $u_*\in\mathcal{E}$ is the number of positive eigenvalues of the linearized operator at such an equilibrium. 
Both the zero number and Morse index can be computed from a permutation of the equilibria, as it was done in \cite{FuscoRocha} and \cite{FiedlerRocha96} for the semilinear dissipative case, see \cite{LappicyQuasi}. For the unbounded structure, a permutation can be computed as in Pimentel and Rocha \cite{RochaPimentel16}. This permutation is called the \emph{Sturm Permutation}. 

We say that two different equilibria $u_-\in \mathcal{E}^b$ and $u_+\in\mathcal{E}=\mathcal{E}^b\cup\mathcal{E}^\infty$ of \eqref{PDE} are \emph{adjacent} (see \cite{Wolfrum02}) if there does not exist any equilibrium $u_*\in\mathcal{E}^b$ of \eqref{PDE} such that $u_*(0)$ lies between $u_-(0)$ and $u_+(0)$, and
\begin{equation}\label{adjacentzeronumb}
    z(u_--u_*) = z(u_--u_+) = z(u_+-u_*).
\end{equation}
Note that if $u_+\in \mathcal{E}^\infty$, the symbol $u_+$ in \eqref{adjacentzeronumb} stands for the $\chi$-component in \eqref{CptEq}.
%
\begin{thm}\emph{\textbf{Unbounded attractor with Chafee-Infante network at infinity.}} \label{attractorthm}
Let $a,f\in C^2$ be bounded such that $a$ is strictly parabolic with limiting behavior \eqref{alim}. Fix $b>0$ such that $\sqrt{b/a^\infty}\not\in \mathbb{N}$. Suppose bounded equilibria are hyperbolic.
Then, 
\begin{enumerate}
    \item There is a bounded heteroclinic orbit, $u(t)\in \mathcal{H}$, between two bounded equilibria, $e_j,e_k\in\mathcal{E}$, i.e.,
    \begin{equation}
        e_j\xleftarrow{t\to -\infty}u(t)\xrightarrow{t\to \infty} e_k,
    \end{equation}
    if, and only if, $e_j$ and $e_k$ are adjacent and $i(e_j)>i(e_k)$.
    
    \item There is a grow-up solution $u(t)$ such that $\mathcal{P}(u(t))\in \mathcal{H}^{up}$ is a heteroclinic from the projected equilibria $\mathcal{P}(e_j)\in\mathcal{P}(\mathcal{E})$ to the equilibria at infinity $\Phi_{k}\in\mathcal{E}^\infty$, i.e.,
    \begin{equation}\label{Hetup}
        \mathcal{P}(e_j)\xleftarrow{t\to -\infty} \mathcal{P}(u(t))\xrightarrow{t\to \infty} \pm\Phi_{k}\footnote{In coordinates, $\mathcal{P}(u(t))\xrightarrow{t\to \infty} \pm\Phi_{k}$ translates to $\left(\frac{u}{ \sqrt{1+||u||_{L^2}^2}},\frac{1}{ \sqrt{1+||u||_{L^2}^2}}\right)\xrightarrow{t\to \infty} (\pm\varphi_k,0)$ in $C^1$.},
    \end{equation}
    if, and only if, $e_j$ and $\Phi_{k}$ are adjacent.
    
    \item There is a heteroclinic at infinity, $\Phi(t)\in \mathcal{H}^\infty\subseteq \mathcal{S}^\infty$, between two equilibria at infinity, $\Phi_{j},\Phi_{k}\in\mathcal{E}^\infty$, i.e.,
    \begin{equation}\label{HetsInfty}
        \Phi_{j}\xleftarrow{t\to -\infty} \Phi(t)\xrightarrow{t\to \infty} \Phi_{k},
    \end{equation}
    if, and only if, $j>k$.
\end{enumerate}
\end{thm}
%
If $b<a^\infty$, i.e. $  N^\infty:= \lfloor\sqrt{b/a^\infty}\rfloor=0$, there are two equilibria at infinity, $\mathcal{E}^\infty=\{\pm \Phi_0\}$. 
If $b>a^\infty$, then the number of equilibria at infinity, $N^\infty$, increases as $b$ increases. Thus the dimension of the attainable subset within the sphere at infinity $\mathcal{S}^\infty$ also increases. 
See Figure \ref{FIG:CI}.
This indicates that grow-up in degenerate diffusion equations (i.e. as $a^\infty$ decreases) may amount to an infinite dimensional attainable sphere at infinity.
%
\begin{figure}[H]
\vspace{-0.25cm}
\minipage[b]{0.33\textwidth}\centering
\begin{subfigure}\centering
    \begin{tikzpicture}[scale=0.75]
    \draw[gray,-] (-2,0) -- (2,0);
    \draw[gray,->] (-1,0) -- (-1.01,0);
    \draw[gray,->] (1,0) -- (1.01,0);
    
    \filldraw [very thick] (-2,0) circle (2pt) node[left] {\footnotesize{$-\Phi_0$}};
    \filldraw [very thick] (2,0) circle (2pt) node[right] {\footnotesize{$+\Phi_0$}};
    \filldraw [gray,very thick] (0,0) circle (1pt) node[above] {\footnotesize{$\mathcal{A}^b$}};

    \filldraw [white,rotate=90,very thick] (-2,0) circle (2pt) node[below] {\footnotesize{$-\Phi_1$}};
    \filldraw [white,rotate=90,very thick] (2,0) circle (2pt) node[above] {\footnotesize{$+\Phi_1$}};   

    \end{tikzpicture}
    \addtocounter{subfigure}{-1}\captionof{subfigure}{\footnotesize{$N^\infty=0$.}}\label{FIG:dim1}
\end{subfigure}
\endminipage\hfill
\minipage[b]{0.33\textwidth}\centering

\begin{subfigure}\centering
    \begin{tikzpicture}[scale=0.75]
    \draw[gray,-] (-2,0) -- (2,0);
    \draw[gray,->] (-1,0) -- (-1.01,0);
    \draw[gray,->] (1,0) -- (1.01,0);
    
    \draw[gray,rotate=90, -] (-2,0) -- (2,0);
    \draw[gray,rotate=90, ->] (-1,0) -- (-1.01,0);
    \draw[gray,rotate=90, ->] (1,0) -- (1.01,0);
    
    \filldraw [rotate=90,very thick] (-2,0) circle (2pt) node[below] {\footnotesize{$-\Phi_1$}};
    \filldraw [rotate=90,very thick] (2,0) circle (2pt) node[above] {\footnotesize{$+\Phi_1$}};    
    
    \draw [very thick] (0,0) circle (57pt);

    \draw[very thick,<-] (1.2,1.6) arc (44:45:0.4cm and 0.4cm);    
    \draw[rotate=180,very thick,<-] (1.2,1.6) arc (44:45:0.4cm and 0.4cm);    

    \draw[rotate=70,very thick,->] (1.2,1.6) arc (44:45:0.4cm and 0.4cm);    
    \draw[rotate=70+180,very thick,->] (1.2,1.6) arc (44:45:0.4cm and 0.4cm);    

    \draw [gray,domain=0:6.28,variable=\t,smooth] plot ({-2+0.64*\t},{0.9*(sin(\t r))});
    \draw[gray,->] (-1,0.9) -- (-1.01,0.9);
    \draw[gray,->] (1,-0.9) -- (1.01,-0.9);
    
    \draw [gray,domain=0:6.28,variable=\t,smooth] plot ({-2+0.64*\t},{-0.9*(sin(\t r))});
    \draw[gray,->] (-1,-0.9) -- (-1.01,-0.9);
    \draw[gray,->] (1,0.9) -- (1.01,0.9);
    
    \filldraw [gray,very thick] (0,0) circle (1pt) node[anchor=south]{\footnotesize{$\mathcal{A}^b$}};
    \filldraw [very thick] (-2,0) circle (2pt) node[left] {\footnotesize{$-\Phi_0$}};
    \filldraw [very thick] (2,0) circle (2pt) node[right] {\footnotesize{$+\Phi_0$}};
    \end{tikzpicture}
    \addtocounter{subfigure}{-1}\captionof{subfigure}{\footnotesize{$N^\infty=1$.}}\label{FIG:dim2}
\end{subfigure}
\endminipage\hfill
\minipage[b]{0.33\textwidth}\centering
\begin{subfigure}\centering
    \begin{tikzpicture}[scale=0.75]
    \draw[gray,dashed,-] (-2,0) -- (2,0);
    \draw[gray,dashed,->] (-1,0) -- (-1.01,0);
    \draw[gray,dashed,->] (1,0) -- (1.01,0);
    
    \draw[gray,rotate=90, dashed,-] (-2,0) -- (2,0);
    \draw[gray,rotate=90, dashed,->] (-1,0) -- (-1.01,0);
    \draw[gray,rotate=90, dashed,->] (1,0) -- (1.01,0);

    \draw[gray,dashed,-] (-0.4,-0.4) -- (0.4,0.4);
    \draw[gray,->] (0.3,0.3) -- (0.31,0.31);
    \draw[gray,rotate=180,->] (0.3,0.3) -- (0.31,0.31);
    
    \filldraw [rotate=90,very thick] (-2,0) circle (2pt) node[below] {\footnotesize{$-\Phi_1$}};
    \filldraw [rotate=90,very thick] (2,0) circle (2pt) node[above] {\footnotesize{$+\Phi_1$}}; 
    \filldraw [gray,very thick] (0,0) circle (1pt);
    \filldraw [very thick] (-2,0) circle (2pt) node[left] {\footnotesize{$-\Phi_0$}};
    \filldraw [very thick] (2,0) circle (2pt) node[right] {\footnotesize{$+\Phi_0$}};   
    \filldraw [very thick] (-0.4,-0.4) circle (2pt) node[anchor=north east] {\footnotesize{$+\Phi_2$}};
    \filldraw [very thick] (0.4,0.4) circle (2pt) node[anchor=south west] {\footnotesize{$-\Phi_2$}};        
    
    \draw [very thick] (0,0) circle (57pt);

    \draw[very thick,<-] (1.2,1.6) arc (44:45:0.4cm and 0.4cm);    
    \draw[rotate=180,very thick,<-] (1.2,1.6) arc (44:45:0.4cm and 0.4cm);    

    \draw[rotate=70,very thick,->] (1.2,1.6) arc (44:45:0.4cm and 0.4cm);    
    \draw[rotate=70+180,very thick,->] (1.2,1.6) arc (44:45:0.4cm and 0.4cm);       

    \draw[very thick] (-2,0) arc (180:360:2cm and 0.4cm);
    \draw[very thick,dashed] (-2,0) arc (180:0:2cm and 0.4cm);

    \draw[rotate=-90, very thick] (-2,0) arc (180:360:2cm and 0.4cm);
    \draw[rotate=-90, very thick,dashed] (-2,0) arc (180:0:2cm and 0.4cm);
    
    \draw[very thick,->] (0.34,1) arc (44:45:0.09cm and 0.4cm);    
    \draw[rotate=90,very thick,->] (0.34,1) arc (44:45:0.09cm and 0.4cm);    
    \draw[rotate=180,very thick,->] (0.34,1) arc (44:45:0.09cm and 0.4cm);    
    \draw[rotate=270,very thick,->] (0.34,1) arc (44:45:0.09cm and 0.4cm);    

    \draw[very thick,->] (-0.34,1) arc (131:130:0.09cm and 0.4cm);    
    \draw[rotate=90,very thick,->] (-0.34,1) arc (131:130:0.09cm and 0.4cm);    
    \draw[rotate=180,very thick,->] (-0.34,1) arc (131:130:0.09cm and 0.4cm);    
    \draw[rotate=270,very thick,->] (-0.34,1) arc (131:130:0.09cm and 0.4cm);    
    \end{tikzpicture}
    \addtocounter{subfigure}{-1}\captionof{subfigure}{\footnotesize{$N^\infty=2$.}}\label{FIG:dim3}
\end{subfigure}
\endminipage
\captionof{figure}{The extended attractor $\overline{\mathcal{P}(\mathcal{A})}\subseteq \mathbb{S}_+$. The projected bounded solutions, $\mathcal{A}^b:=\mathcal{E}^b\cup \mathcal{H}^b$, are represented by a (gray) point. The projected grow-up solutions (gray arrows), $\mathcal{H}^{up}$, are heteroclinics from projected bounded solutions towards the sphere at infinity (black), $\mathcal{S}^\infty$. The attainable $N^\infty$-dimensional sphere within $\mathcal{S}^\infty$ has a Chafee-Infante type structure with equilibria at infinity $\Phi_j$ and their heteroclinics.
}\label{FIG:CI}
\end{figure}
The remaining is organized as follows. Section \ref{sec:func} provides the necessary background for the main proof, which is provided in Section \ref{sec:pf}. 
We compute the bounded attractor in Subsection \ref{sec:bddsturm}, 
we describe the induced flow at infinity in Subsection \ref{sec:unbddsturmlin}, 
and we describe the grow-up solutions in Subsection \ref{sec:growup}. 
Lastly, we include concluding remarks in Section \ref{sec:disc}. 

\section{Background} \label{sec:func}

This section provides all the necessary tools for the proof of the main theorems. 
First, we use semigroup theory to guarantee not only that solutions exist, but also that some of them grow up. 
Second, we present the nodal property (also known as dropping lemma), which guarantees that the number of intersections of the spatial profiles of solutions of \eqref{PDE} do not increase in time.
Lastly, we introduce the Fiedler-Brunovský-map (a.k.a. $y$-map), and its major consequence: prescribe dropping times (when the number of intersections of solutions drops) with realizable initial data that achieves a dropping of intersections at such prescribed times.

\subsection*{Semigroup theory}

The phase-space $X^\alpha$ lies in the space of H\"older continuous functions $X:=C^\beta([0,\pi])$ with H\"older coefficient $\beta\in (0,1)$, intersected with the Neumann boundary conditions, constructed as follows. See \cite{Lunardi95,Amann95,BabinVishik92}. The notation $C^{\beta}$ for some $\beta\in\mathbb{R}_+$ indicates that $\beta$ can be rewritten as $\lfloor\beta\rfloor+\{\beta\}$, where the integer part $\lfloor\beta\rfloor\in\mathbb{N}$ denotes that functions are up to $\lfloor\beta\rfloor$-times differentiable, and its $\lfloor\beta\rfloor$-derivative is $\{\beta\}$-Hölder, where $\{\beta\}\in[0,1)$ is the fractional part of $\beta$.

Equation \eqref{PDE} can be rewritten as the following abstract equation,
\begin{equation}\label{linequiv}
    u_t=Au+bu+g(u)    
\end{equation}
where $A:D(A)\rightarrow X$ is the linearization of the right-hand side of \eqref{PDE} with $b=0$, at any point in the neighborhood of the initial data $u_0(x)$, and $g$ is the Nemitskii operator of the remaining terms, which takes values in $X$, given by $g(u):=a(x,u,u_x)u_{xx}+f(x,u,u_x)-Au$. 
The domain of $A$ is $D(A):=C^{2,\beta}([0,\pi])\cap \{\text{Neumann b.c.}\}\subseteq X$, where $\beta\in (0,1)$. Moreover, consider the interpolation spaces $X^\alpha:=C^{2\alpha+\beta}([0,\pi])$ between $X$ and $D(A)$, with $\alpha\in (0,1)$, such that $A$ generates a strongly continuous semigroup in $X^\alpha$. Therefore, solutions of the equation \eqref{PDE} define a semiflow in $X^\alpha$ according to the variation of constants formula. In particular, this settles existence and uniqueness matters. See Lunardi \cite[Theorem 8.1.1]{Lunardi95}.

We suppose $2\alpha+\beta>1$ so that solutions are in $C^1([0,\pi])$. Moreover, $C^{2\alpha+\beta}\subseteq L^2$, and hence the 
subspace $C^{2\alpha+\beta}$ inherits the inner product $\langle\cdot,\cdot\rangle$ of $L^2$. From now on, all norms and inner products are the usual ones in $L^2$, except when explicitly expressed.
Note that the operator $\partial^2_x$ has eigenvalues $\lambda_j=-j^2$ , with $j=0,1,\dots$, and eigenfunctions $\varphi_j(x)=\cos(jx)$, which form an orthonormal basis of $L^2$, and thereby of $X^\alpha$. Therefore, we can decompose the dynamics by a semiflow in each eigendirection. In the next Lemma, we guarantee that solutions exist for all time, yet, they become unbounded in infinite time.

\begin{lem}\label{lemnondiss}
If $b>0$, then solutions $u(t)$ of \eqref{PDE} are global in time. Moreover, there exists solutions $u(t)$ that grow up (blow up in infinte time) in $X^\alpha$. 
\end{lem}
\begin{proof}
Firstly, note that since $a,f\in C^2$ are bounded with uniformly parabolic $a\geq \epsilon$, the semigroup $u(t)$ given by the variation of constants formula is bounded for any finite time. Therefore, we can indefinitely extend any solution in time, and the maximal time of existence is $T=\infty$. Hence, finite time blow-up does not occur. 

Next, we show that some solutions indeed become unbounded in $X^\alpha$ as $t\to\infty$. 
Consider sub-solutions $\underbar{$u$}(t,x)$ of \eqref{PDE} that satisfy the auxiliary semilinear problem, due to the strict parabolicity condition $a(x,u,u_x)\geq \epsilon >0$,
\begin{equation}\label{subsol}
    \underbar{$u$}_t=\epsilon \underbar{$u$}_{xx}+b\underbar{$u$}+f(x,\underbar{$u$},\underbar{$u$}_x).
\end{equation}
with the same initial data as \eqref{PDE}, $\underbar{$u$}_0(x)=u_0(x)$, and Neumann boundary conditions. Hence, $u$ and $\underbar{$u$}$ coincide in the parabolic boundary, and the comparison principle implies that $u(t,x)\geq \underbar{$u$}(t,x)$ for all $(t,x)\in\mathbb{R}_+\times [0,1]$, see \cite[Theorem A.10]{ArrietaCarvalhoBernal00}. 
Note that to apply such comparison, we can write the quasilinear equation \eqref{PDE} abstractly as \eqref{linequiv} and compare $g(u)\geq g_\epsilon (u):=\epsilon u_{xx}+f(x,u,u_x)-Au$. 

Moreover, due to \cite[Lemma 1]{RochaPimentel16}, there is a grow-up (in $L^2$) solution $\underbar{$u$}(t,x)$ for $b>0$, then $u(t,x)$ must also grow up (in $L^2$). 
Lastly, the embedding $X^\alpha \subseteq L^2$ implies that $||.||_{L^2} \leq c ||.||_{\alpha}$ for some $c\in\mathbb{R}$, which in turn implies that $u(t)$ also grows up in $X^\alpha$.
\end{proof}
Therefore, we have proved that solutions in $X^\alpha$ exist globally in time, even though some solutions may grow up. The projection $\mathcal{P}$ in \eqref{P} and its inverse, $\mathcal{P}^{-1}(\chi,z):=(\chi / z,1)$, thereby induce a semiflow in the relevant subset within the Poincaré hemisphere, $\mathcal{P}(X^\alpha\times \{1\})\subseteq \mathbb{S}_+|_{z>0}$, which evolves according to the equation \eqref{flowSPHERE} for $z>0$.
Next, we show in Lemma~\ref{lem:EQUADORinv} that such induced semiflow is well-defined in the limit $z\to 0$ and the sphere at infinity $\mathcal{S}^\infty:=\mathbb{S}_+|_{z=0}$ is invariant. 

\begin{lem}\label{lem:EQUADORinv}
The semiflow of equation \eqref{flowSPHERE} is well-defined in the limit $z\to 0$ and the sphere at infinity $\mathcal{S}^\infty$ is invariant.
\end{lem}
\begin{proof}
In order to prove that one can pass to the limit in equation \eqref{flowSPHERE}, one has to guarantee that $\langle \mathcal{L}^z(\chi),\chi \rangle$ remains bounded, as $z\rightarrow 0$, yielding a well-defined vector field at $z=0$. 
Indeed,
\begin{equation}\label{AEAE}
    |\langle \mathcal{L}^z(\chi),\chi \rangle |  
    \leq   a||\chi_{x}||^2+b||\chi||^2+\langle f^z(x,\chi,\chi_x),\chi \rangle,
\end{equation}
since the quasilinear diffusion coefficient is bounded by some $a\in\mathbb{R}_+$. 

To guarantee that \eqref{AEAE} is bounded, we only need to prove that $\chi$ and $\chi_x$ have bounded $L^2$ norms.
By definition of $\chi$ in \eqref{P}, it follows that $\chi \in\mathbb{S}_+$ and hence $\Vert \chi \Vert=1$. In order to prove that $\Vert\chi_x\Vert$ is bounded, it is enough to guarantee that $\Vert\chi\Vert_\alpha$ is bounded, since: $\Vert \chi_x \Vert \leq c\Vert \chi_x \Vert_{X^\alpha}$, due to the embedding $X^\alpha\subseteq L^2$.

Next, we proceed to prove that $\Vert\chi\Vert_\alpha$ is bounded.
First, we decompose $u(t)$ in the basis $\{  \varphi_{j} \}_{j\in\mathbb{N}}$, in order to apply \cite[Lemma 5.1]{P}. More precisely, we denote by $p(t)$ the sum of the (finitely many) terms related to growth, and by $q(t)$ the sum of remaining terms related to decay (which remain uniformly bounded), as time goes to infinity. Thus
\begin{equation}\label{ultimodia}
        u(t)=p(t)+q(t):=\sum_{j=0}^{M} u_j(t) \phi_{j}+\sum_{j>M} u_j(t) \phi_{j},
\end{equation}
for some integer $M$. Therefore, 
\begin{equation}
    \chi(t)=\frac{p(t)}{\sqrt{1+\Vert u\Vert^2}}+\frac{q(t)}{\sqrt{1+\Vert u\Vert^2}}.
\end{equation}
Since the second term goes to zero as $t\to \infty$, and the first term lies in a finite dimensional space, where the norms $\Vert .\Vert$ and $\Vert .\Vert_\alpha$ are equivalent, we obtain
\begin{equation}
    \Vert\chi\Vert_\alpha\leq  C \left\Vert\frac{p(t)}{\sqrt{1+\Vert u\Vert^2}}\right\Vert_\alpha\leq \tilde{C} \left\Vert\frac{p(t)}{\sqrt{1+\Vert u\Vert^2}}\right\Vert <\infty    
\end{equation}
where the last norm is bounded since $\Vert \chi \Vert$ is bounded. 
This finishes the proof that $\langle \mathcal{L}^z(\chi),\chi \rangle$ in \eqref{AEAE} is bounded.

We can now compute the limit of the equation \eqref{flowSPHERE} as $z\to 0$, which is given by
\begin{subequations}\label{z=0EEpf}
\begin{align}
    \chi_t&=a^\infty\chi_{xx}-\langle a^\infty\chi_{xx},\chi\rangle \chi,\\
    z_t&=0
\end{align}
\end{subequations}
due to \eqref{alim}. Note that $\chi$ has Neumann boundary conditions. 

In order to guarantee the existence of a semiflow for \eqref{z=0EEpf}, note that the linear term, $a^\infty\chi_{xx}$, generates an analytic semigroup, whereas nonlinearity can be rewritten as $a^\infty||\chi_{x}||^2 \chi$, which is the product of two Lipschitz functions.
Moreover, equation \eqref{z=0EEpf} implies invariance of $\mathcal{S}^\infty$, as $z_t=0$.
\end{proof}

\subsection*{Nodal Property}

Let the \emph{zero number} $0\leq z(u(t,.))\leq \infty$ count the number of strict sign changes of the spatial (i.e., with fixed $t$) profiles $u(t,.):[0,\pi]\to\mathbb{R}\cup\{\pm\infty\}$. More precisely, if $x\mapsto u(t,x)$ is not of constant sign, let
\begin{equation} \label{zt}
    z(u(t,.)):= \sup_k \left\{  
        \begin{array}{c} 
        \text{There is a partition $\{ x_j\}_{j=1}^{k}$ of } [0,\pi]\\
        \text{such that } u(t,x_j)u(t,x_{j+1})<0 \text{ for all } j=0,\dots,k-1
        \end{array} \right\}.
\end{equation}
For functions which do not change sign, $x\mapsto u(t,x)\neq 0$, we define $z(u):=0$. For the trivial constant, $x\mapsto u(t,x)\equiv 0$, we define $z(u\equiv 0):=-1$. 
Note we allow discontinuous and unbounded profiles $u$.

Now we present the nodal property. Different versions of this well known fact are due to Sturm \cite{Sturm}, Matano \cite{Matano82}, Angenent \cite{Angenent88} and others.
We recall that a point $(t_0,x_0)\in\mathbb{R}\times [0,\pi]$ such that $u(t_0,x_0)=0$ is said to be a \emph{simple zero} of $u(t,.)\in C^1$ if $u_x(t_0,x_0)\neq 0$ and a \emph{multiple zero} if $u_x(t_0,x_0)=0$.

\begin{lem} \label{droplem}
    \emph{\textbf{Dropping lemma.}}
    Consider a non-trivial $C^1$ solution of the equation 
    \begin{equation}\label{linPDE}
        v_t=a(t,x)v_{xx}+b(t,x)v_x+c(t,x)v
    \end{equation}
    where $x\in [0,\pi]$ has Neumann boundary conditions, $a,b,c\in L^\infty$ for $t\in [0,T)$ and $a(t,x)>0$ for all $(t,x)\in [0,T)\times [0,1]$. Then, the zero number $z(v(t,.))$ satisfies
    \begin{enumerate}
    \item $z(v(t,.))<\infty$ for any $t\in (0,T)$;
    \item $z(v(t,.))$ is nonincreasing with time $t$;
    \item $z(v(t,.))$ decreases at multiple zeros $(t_0,x_0)$ of $v(t,.)$, i.e., 
    \begin{equation}
        z(v(t_0-\epsilon))>z(v({t_0+\epsilon}))
    \end{equation} 
    for any  sufficiently small $\epsilon>0$.
    \end{enumerate}
\end{lem}
This result proves that the number of intersections of two different solutions of \eqref{PDE} is nonincreasing with time $t$, and decreases whenever a multiple zero occurs.
Indeed, the difference $v:=u^1-u^2$ of any two solutions $u_1,u_2$ of the nonlinear equation \eqref{PDE} satisfies a linear equation of the type \eqref{linPDE}, where
\begin{subequations}
\begin{align}
    {a}(t,x)&:=\int_0^1 a(x,u^s,u_x^s )ds,\\
    {b}(t,x)&:=\int_0^1 f_p(x,u^s ,u_x^s )ds,\\
    {c}(t,x)&:=\int_0^1 f_u(x,u^s ,u_x^s )ds,   
\end{align}
\end{subequations}
where $p:=u_x$, and $u^s :=(1-s) u_1 +s u_2 $ for $s\in[0,1]$. 

\subsection*{Fiedler-Brunovsk\'y-map}
Now we introduce the \emph{y-map}, 
\begin{equation}
    y:D(y)\subseteq X^\alpha \rightarrow \mathbb{S}^n,
\end{equation}
a tool used to detect heteroclinics for dissipative equations firstly used by Brunovsk\'y and Fiedler in \cite{FiedlerBrunovsky88}, and later on for non-dissipative equations possessing grow-up in \cite{BenGal10,RochaPimentel16}.
Roughly speaking, $y(u_0)$ encodes information about the zero number $z(u(t,.))$ of a solution $u(t,.)$ of \eqref{PDE} with initial data $u_0\in X^\alpha$ for all times $t\in\mathbb{R}_+$. We will show in this section that surjectivity of $y$ for suitable domains $D(y)$ implies that there is an initial data $u_0\in D(y)$ 
with prescribed numbers of zeroes and dropping times.

Consider a solution $v(t)$ of \eqref{linPDE} with initial data $v_0$ such that $z(v_0)=n$. Hence, $z(v(t,.))$ can drop at most $n$ times, due to the Lemma \ref{droplem}. For each $k=0,...,n$, define the \emph{dropping times} $t_k\in\mathbb{R}_+$ as the first time that the zero number of the solution $v(t)$ drops to $k$ zeros or less, i.e.,
\begin{equation}
    t_k:=\inf \{t\in\mathbb{R}_+: z(v(t,.))\leq k\}.
\end{equation}
In case the number of zeros of the solution $v(t)$ does not pass below the $k$-level, i.e., if $z(v(t,.))>k$ for any $t\in\mathbb{R}_+$, we say that the dropping occurs at infinite time, and denote by $t_k:=\infty$.
Next, we compactify the dropping times $t_k$ via 
\begin{equation}
    \tau_k:=\tanh(t_k)\in [0,1].
\end{equation}
Note $\tau_n=t_n=0$. 
Moreover, $0=\tau_n\leq ...\leq \tau_k \leq ... \leq \tau_{0}$, 
due to the Lemma \ref{droplem}. Note that the inequalities are not strict, since the dropping does not necessarily occur one-by-one: the zero number can drop by more than one at a single dropping time.

Note that $v(t,0)\neq 0$ for all $t\in (t_k,t_{k-1})$, if $t_k\neq t_{k-1}$.\footnote{Indeed, if this was the case, and $v(t_*,0)=0$ for some $t_*\in (t_k,t_{k-1})$, then $(t_*,0)$ would be a multiple zero, due to Neumann boundary conditions. The Lemma \ref{droplem} would imply that $t_*$ is a dropping time, which is not in the list $t_k$, yielding a contradiction.} Hence, the sign of $v(t,0)$ is constant for $t\in (t_k,t_{k-1})$, since $v(t,x)$ is $C^0$ in $t$ and $C^1$ in $x$. We then define the sign of the constant boundary values between two dropping times for some $t_*\in (t_k,t_{k-1})$:
\begin{equation}\label{iota}
    \iota_k:=
    \begin{cases}
        \text{sign}(v(t_*,0)),  & \text{ if } t_k<t_{k-1},\\
        0, & \text{ if } t_k=t_{k-1}.
    \end{cases}
\end{equation}
%

We now define the $y$-map, and then discuss its image and domain $D(y)$, respectively.
For any initial data which has at most $n$ zeros, $u_0\in \{ u_0\in X^\alpha : z(u_0)\leq n\}\backslash \{0\}$, we define the coordinate $y_k$ of the $y$-map $y=(y_{0},...,y_n)\in \mathbb{R}^{n+1}$ by
\begin{equation}\label{ymapcomp}
    y_k(u_0):=\iota_k \sqrt{\tau_{k-1}-\tau_k}
\end{equation}
for each $k=0,...,n$. Also, $\tau_{-1}:=1$ is chosen for well-definition purposes. Continuity of $y$ with respect to $f,u_0$ and $t$ was proved in \cite[Lemma 2.1]{FiedlerBrunovsky88}.
Similar arguments can be replicated for equation \eqref{PDE}, including continuity in $a$.

The image of the $y$-map lies in the $n$-sphere, $\mathbb{S}^n\subseteq \mathbb{R}^{n+1}$, since 
\begin{equation}
    ||y||_{\mathbb{R}^{n+1}}=\sum_{k=0}^n y_k^2=\sum_{k=0}^n (\tau_{k-1}-\tau_k)=\tau_{-1}-\tau_n=1.
\end{equation}
To obtain a map between $n$-dimensional spheres, we restrict the domain $D(y)$ as follows. 
Consider 
the linearization of \eqref{PDE} at a hyperbolic equilibrium $u_*$ with Morse index $i(u_*)=n$, which yields a linear operator $A_*$ with simple eigenvalues $\{ \lambda^*_{k}\}_{k\in\mathbb{N}_0}$ that acumulate at $-\infty$, and corresponding eigenfunctions $ \{ \varphi^*_{k} \}_{k\in \mathbb{N}_0}$ such that $z(\varphi^*_{k})=k$, due to Sturm-Liouville theory. 
Hence, there is an $n$-dimensional unstable manifold, $W^u(u_*)\subseteq X^\alpha$, which is a graph over its tangent space, given by $E^u:=\text{span} \{\varphi_{k}\}_{k=0}^{n-1}$, via a diffeomorphism $h:E^u\to W^u(u_*)$.
We consider a sphere centered at the origin with sufficiently small $\delta>0$ radius, $\mathbb{S}^n_\delta\subseteq E^u$, which is lifted by $h$ to an 
$n$-sphere within the unstable manifold of $u_*$, denoted by $\Sigma^n:=h(\mathbb{S}^n_\delta)$. The domain of the $y$-map is $D(y):=\Sigma^n$. In particular, note that $z(u_0)\leq n$ for any $u_0\in \Sigma^n$, due to \cite{BrunovskyFiedler86}.

The image $y(u_0)\in \mathbb{S}^n$ encodes the information about the zero numbers $z(u(t,.))$ for a solution $u(t,.)$ of \eqref{PDE} with initial data $u_0\in \Sigma^n$ for all $t\in\mathbb{R}_+$. 
To prove surjectivity of the $y$-map, we show it is \emph{essential}, i.e., it is not homotopic to a constant map.
\begin{lem} \label{ymaplem}
    The map $y:\Sigma^n\to \mathbb{S}^n$ with coordinates given by \eqref{ymapcomp} is essential. In particular, it is surjective.
\end{lem}
\begin{proof}
Essentiality implies surjectivity. Else, for some $y_*\in \mathbb{S}^n$, there would not exist any $u_0\in D(y)$ such that $y_*=y(u_0)$. Therefore $y(\Sigma^n)\subseteq \mathbb{S}^n\backslash \{y_*\}$. Note $\mathbb{S}^n\backslash \{y_*\}$ is homeomorphic to a ball of dimension $n$, which is contractible to a point. Hence the image of the $y$-map would be homotopic to a constant map, a contradiction.

Note that essentiality is preserved throughout homotopies. Hence, to prove that the $y$-map is essential with respect to nonlinear semiflows $u(t)$ generated by \eqref{PDE}, we homotope the nonlinear semiflow to a linearized semiflow $v(t)$ generated by \eqref{linPDE}, and prove that the $y$-map is essential with respect to linear semiflows.
Indeed, consider
\begin{equation}
    u_t=a^\tau(x,u,u_x)u_{xx}+bu+f^\tau(x,u,u_x)
\end{equation}
where the homotopy diffusion is $a^\tau(x,u,u_x):=\tau a(x,u,u_x)+(1-\tau)a(x,0, 0)$, whereas the homotopy reaction is $f^\tau(x,u,u_x):=\tau f(x,u,u_x)+(1-\tau)[f_p(x,0, 0)u_x+f_u(x,0, 0)u]$,
for $\tau\in[0,1]$.
Note the continuity of all terms with respect to $\tau$.
Moreover, the linearized equation at $0$ remains unchanged throughout the homotopy. 
Hence the uniform hyperbolicity of $0$ is guaranteed during the homotopy and it has the same Morse index for all homotopy parameter $\tau\in[0,1]$.

Lastly, note that the $y$-map is odd for linear semiflows, and consequently essential by the Borsuk-Ulam Theorem. 
\end{proof}

\begin{cor}\label{ymapcor}
    Let $u_*$ be a hyperbolic equilibrium of \eqref{PDE} with Morse index $i(u_*)$. Consider the set $\Sigma^n \subseteq W^u(u_*)\backslash \{u_*\}$ which is homotopic to an $n$-dimensional sphere, with $n:=i(u_*)-1$, centered at $u_*$.
    
    Then, for any sequences $0=t_n\leq ...\leq t_0\leq \infty$ and $\iota_n,...,\iota_0\in \{\pm 1\}$, there is an initial datum $u_0 \in \Sigma^n$ with corresponding solution $u(t)$ of \eqref{PDE} such that the graph $t\mapsto z(u(t)-u_*)$ is determined by $\{t_k\}_{k=0}^n$ as follows
    \begin{subequations}
    \begin{align}
        z(u(t)-u_*)&\leq k \text{ for all $t\geq t_k$}\\
        \mathrm{sign}(u(t)-u_*)&=\iota_k \text{ at $x=0$, for all $t\in (t_k,t_{k-1})$}.
    \end{align}
    \end{subequations}
\end{cor}

\section{Proof}\label{sec:pf}

\subsection{Bounded Sturm structure}\label{sec:bddsturm}

In this section, we cut off the quasilinear parabolic equation \eqref{PDE} in order to obtain a characterization of the maximal compact attractor within the unbounded attractor.

Define $\mathcal{A}^{cpt}$ to be the maximal compact invariant set within $\mathcal{A}$. Therefore, it consists of the global bounded solutions, and its complement is the unbounded part of the attractor, 
which consists of the unbounded solutions. 
Let $c>0$ denote a bound for the solutions $u\in \mathcal{A}^{cpt}\subset X^\alpha\subset C^1$, i.e., 
\begin{equation}
       \sup_{t\in\mathbb{R},x\in [0,\pi], u\in \mathcal{A}^{cpt}}\{ |u(t,x)|,|u_x(t,x)| \}<c. 
\end{equation}
Consider the set $N_{\mathcal{A}^{cpt}}:=\{ (x,u,u_x)\in [0,\pi]\times\mathbb{R}^2: |u|,|u_x|\leq c  \}$ and an open neighborhood $N_\epsilon\subset [0,\pi]\times\mathbb{R}^2$ of $N_{\mathcal{A}^{cpt}}$. We modify $f$ outside $N_{\mathcal{A}^{cpt}}$ as follows
\begin{equation}
F(x,u,u_x):=
\begin{cases}
bu+f(x,u,u_x)  & \text{ for } (x,u,u_x)\in N_{\mathcal{A}^{cpt}}\\
G(x,u,u_x) & \text{ for } (x,u,u_x)\in N_\epsilon\setminus N_{\mathcal{A}^{cpt}}\\
-u & \text{ for } (x,u,u_x)\notin N_\epsilon,
\end{cases}
\end{equation}
where $G$ is a transition function obtained from Urysohn's lemma. 

Therefore, solutions of \eqref{PDE} contained in the set $\mathcal{A}^{cpt}$ coincide with the semiflow of the following modified dissipative equation $u_t=a(x,u,u_x)u_{xx}+F(x,u,u_x)$.
Hence, the characterization of the bounded attractor $\mathcal{A}^{cpt}$ follows the dissipative case: 
it consists of bounded equilibria and bounded heteroclinics between them.
Moreover, the heteroclinics between bounded equilibria of quasilinear equations is constructed as in \cite{LappicyQuasi}. 
This partially proves Theorem \ref{attractorthm1}, and proves item 1 of Theorem \ref{attractorthm}.

\subsection{Sphere at infinity with Chafee-Infante type dynamics}\label{sec:unbddsturmlin} %

In this section, 
we prove that the sphere at infinity, $\mathcal{S}^\infty$, has gradient structure, which completes the proof of Theorem \ref{attractorthm1}. 
Moreover, in Proposition~\ref{lem:unbddeq}, 
we fully construct the structure given by equilibria at infinity, $\mathcal{E}^\infty$, and the heteroclinics at infinity, $\mathcal{H}^\infty$, which proves item 3 of Theorem \ref{attractorthm}. 
%
%

To obtain an accurate description of the (nonlinear and nonlocal) semiflow at the sphere at infinity, $\mathcal{S}^\infty$, we compute the limiting equation \eqref{flowSPHERE1} as $z\to 0$, which can be accomplished due to Lemma \ref{lem:EQUADORinv}. The limit of \eqref{flowSPHERE1} as $z\to 0$ is given by
\begin{equation}\label{flowSPHEREz0}
    \chi_t=a^\infty\chi_{xx}- a^\infty\langle \chi_{xx},\chi\rangle \chi,
\end{equation}
with Neumann boundary conditions, since $a^z(x,\chi,\chi_x)\to a^\infty$ as $z\to 0$, due to \eqref{alim}.
Alternatively, each coordinate $\chi_j=\langle \chi,\varphi_j\rangle$ satisfies the following equation
\begin{equation}\label{flowSPHEREz0coordNEW}
    (\chi_j)_t= a^\infty [\lambda_j+||\chi_{x}||^2 ]\chi_{j}.
\end{equation}
%
%
%
%

\begin{lem}\label{lem:EQUADORLyap}
There is a Lyapunov function for the nonlocal nonlinear semiflow \eqref{flowSPHEREz0} at the sphere at infinity, $\mathcal{S}^\infty$, given by 
$E^\infty:\mathcal{S}^\infty\to \mathbb{R}$, where
\begin{equation}\label{EQUADORLyap}
   E^\infty(\chi)= \frac{||\chi_x||^2}{2} \qquad \text{ such that }  \qquad \frac{dE^\infty}{dt}(\chi(t))\leq 0.
\end{equation}
Moreover, $dE^\infty (\chi(t))/dt=0$ if, and only if, $\chi(t)$ is an equilibrium of \eqref{flowSPHEREz0}.
\end{lem}
\begin{proof}
Differentiating \eqref{EQUADORLyap} with respect to time, along a solution of \eqref{flowSPHEREz0}, integrating by parts, and substituting  the equation \eqref{flowSPHEREz0}, yields 
\begin{subequations}
\begin{align}
    \frac{dE^\infty}{dt}&=- \langle \chi_t,\chi_{xx}\rangle \label{lyap}\\
    &=-a^\infty \left[  || \chi_{xx}||^2-\langle \chi_{xx},\chi \rangle^2 \right] \leq 0.
\end{align}
\end{subequations}
%
where the last inequality holds due to Cauchy-Schwarz, $\langle \chi,\chi_{xx}\rangle \leq ||\chi || \cdot ||\chi_{xx}||$, altogether with the fact that $\chi$ lies on the sphere at infinity, i.e., $||\chi||=1$. 
Moreover, equation \eqref{lyap} implies that critical points of $E^\infty$ correspond to equilibria of \eqref{flowSPHEREz0}.
%
\end{proof}



Next, we consider a secondary projection
(into the tangent spaces of the sphere at infinity), 
so that the induced semiflow is local, and thus the dynamics can be easily described. As a preliminary step, we need to account for how many growing directions an unbounded solution has, and which is the eigendirection with biggest growth rate.
\begin{lem}\label{growingmodes}
Consider a grow-up solution $u(t)$ of \eqref{PDE} with associated eigendirectional semiflows given by $u_j(t):=\langle u(t) , \varphi_j \rangle$ for $j\in\mathbb{N}_0$.
Hence, 
\begin{enumerate}
    \item[(i)] $u_j(t)$ grows for $j=0,\dots,N^\infty$, where $N^\infty= \lfloor \sqrt{b/a^\infty} \rfloor$, in case its initial condition is nonzero, i.e., $u_{j}(0)\neq 0$. Moreover, if two eigenprojections $u_j(t)$ and $u_{j+1}(t)$ such that $j+1\leq N^\infty$ grow up, then $u_{j}(t)$ grows faster than $u_{j+1}(t)$.
    
    \item[(ii)] $u_j(t)$ are bounded for the remaining $j>N^\infty$.
\end{enumerate}
\end{lem} 
\begin{proof}
Since the solution $u(t)$ grows up, it lies outside a ball of $X^\alpha$ with sufficiently large radius $R$, for sufficiently large time. Hence, the Nemitskii operator of $a$ satisfies $a^\infty-\delta \leq a(x,u,u_x)\leq a^\infty+\delta$ for sufficiently small $\delta>0$, due to \eqref{alim}. Therefore, the subsolutions $u^-$ and supersolutions $u^+$ of \eqref{PDE} with respective semilinear diffusion coefficient $a^\infty-\delta$ and $a^\infty+\delta$ satisfy
\begin{equation}\label{OIE}
    (a^\infty-\delta) u^-_{xx}+bu^-+f(x,u^-,u^-_x) \leq u_t \leq (a^\infty+\delta) u^+_{xx}+bu^+ +f(x,u^+,u^+_x).
\end{equation}
We project the semiflow in the $\varphi_j$ eigendirection, which yields
\begin{equation}\label{OIE2}
    ((a^\infty-\delta)\lambda_j+b)u^-_j+f_j(t) \leq (u_j)_t \leq ((a^\infty+\delta)\lambda_j+b)u^+_j+f_j(t),
\end{equation}
where $f_j(t):=\langle f(x,u,u_x),\varphi_j\rangle$. The variation of constants formula leads to
\begin{equation}\label{LB}
       e^{[(a^\infty-\delta)\lambda_j+b]t}u^{-}_j(0) +I^-_j(t)\leq u_j(t) \leq e^{[(a^\infty+\delta)\lambda_j+b]t}u^+_j(0)+I^+_j(t)
\end{equation}
where $u^{\pm}_j(0):=u_j(0)+\int_0^\infty e^{-[(a^\infty\pm\delta)\lambda_{j}+b]s} f_{j}(s)ds$ and $I^\pm_j(t):=\int_\infty^t e^{[(a^\infty\pm\delta)\lambda_{j}+b](t-s)} f_{j}(s)ds$.

Now, we want to collect all indices $j$ such that the subsolutions $u_j^-(t)$ 
grow (in forward time), which occurs when 
\begin{equation}\label{COND1}
    (a^\infty-\delta)\cdot  (-j^2)+b>0.    
\end{equation}
Note that the term $I^-_j(t)$ is bounded, and it does not contribute to the growth of solutions. 
There are finitely many eigendirections so that $j<\sqrt{b/(a^\infty-\delta)}$, for sufficiently small $\delta>0$. Thus, the number of growing eigendirections is $j\leq N^\infty $, which yield growth in case that $u_j(0)\neq 0$. This proves the first part of item $(i)$. 

To compare the growth of $u_{j+1}$ and $u_j$, in case both coordinates grow, note that the growth rates of subsolutions and supersolutions in \eqref{LB} are all different for different $j$. 
Thus, $u_{j}$ grows faster than $u_{j+1}$ if the upper bound of the growth rate of $u_{j+1}$ is smaller than the lower bound of the growth rate of $u_j$, i.e.
%
\begin{equation}
    (a^\infty+\delta)(\lambda_{j+1})<(a^\infty-\delta)(\lambda_{j}).
\end{equation}
Equivalently, $(a^\infty+\delta)(j+1)^2>(a^\infty-\delta)j^2$, which is true for sufficiently small $\delta>0$.
This proves the second part of $(i)$.

Similarly, we collect all indices $j$ such that the supersolutions $u_j^+$ in \eqref{LB} yield a bounded (in forward time) projected eigendirection, which occurs when $(a^\infty+\delta)\cdot  (-j^2)+b<0$. Note that the term $I_j^+$ remains bounded.
Thus, there is an infinite number of bounded eigendirections $u_j(t)$, given by $j>N^\infty$. This proves $(ii)$. 
\end{proof}
Note that a solution $u(t)$ grows up in the $\varphi_j$ eigendirection if, and only if the normalized semiflow in the $j$-eigendirection converges to 1, i.e.,
\begin{equation}\label{growupdirection}
    \lim_{t\to \infty } \frac{u(t)}{\Vert u(t)\Vert}= \pm \varphi_j, \hspace{0.5pt} \text{ in the $L^2$-topology} \iff \lim_{t\to \infty}\frac{u_j(t)}{\Vert u(t)\Vert}=\pm 1.
\end{equation}
The characterization \eqref{growupdirection} follows from the relation:
\begin{equation}
\left\Vert \frac{u(t)}{\Vert u(t) \Vert} - \varphi_j \right\Vert^2= 2\left( 1 - \frac{u_j(t)}{\Vert u(t) \Vert}\right).
\end{equation}
Thus, due to Lemma \ref{growingmodes}, there is an unique $j^*$, given by the smallest index $j\leq N^\infty$ such that $u_{j}(0)\neq 0$, which yields the largest growth rate for a grow-up solutions $u(t)$. Therefore, 
the normalized grow-up solution $u(t)$ converges to $\varphi_j^*$, according to \eqref{growupdirection}. 

Consider a grow-up solution $u(t)$ with biggest growth rate (i.e., nonzero initial data) in the $\varphi_j^*$-eigendirection. Consider also the hyperplane $C_{j^*}$ which is tangent to the sphere at infinity at $(\varphi_{j^*},0)\in\mathcal{S}^\infty$, given by
\begin{equation}
    C_{j^*}:=\{ (\chi,z)\in L^2\times [0,1]\text{ $|$ } \chi_{j^*}=+1,\chi_j\in\mathbb{R}\text{ for all } j\in\mathbb{N}_0 \backslash\{j^*\}\}.
\end{equation}
Similar to the projection $\mathcal{P}$ in \eqref{P}, we will construct another projection $\Tilde{\mathcal{P}}_{j^*}$ from $X^\alpha\times \{1\}$ into $C_{j^*}$. Consider any point $u\in X^{\alpha}\subseteq L^2$ and a line that passes through the points $(u,1),(0,0)\in  L^2\times [0,1]$. See Figure \ref{fig:Ptilde}. The intersection of such line with the plane $C_{j^*}$ defines the projection $\Tilde{\mathcal{P}}_{j^*}$, with coordinates given by $(\xi,\zeta)$:
\begin{equation}\label{Ptilde}
    (\xi,\zeta):=\Tilde{\mathcal{P}}_{j^*}(u,1)=\left(\frac{u}{ \langle u,\varphi_{j^*} \rangle},\frac{1}{ \langle u,\varphi_{j^*} \rangle}\right).
\end{equation}
\begin{figure}[H]\centering
    \begin{tikzpicture}[scale=2]
    \draw[->] (-1.5,0) -- (1.5,0) node[right] {$\xi\in L^2 $};
    \draw[->] (0,-0.1) -- (0,1.5) node[above] {$\zeta\in [0,1]$};

    \draw [thick, domain=0:3.14,variable=\t,smooth] plot ({cos(\t r)},{sin(\t r)}); 
    
    \draw[-] (1.5,1) -- (-1.5,1) node[left] {$L^2\times \{ 1\}$};
    \draw[-] (1,-0.1) -- (1,1.5) node[above] {$C_{j^*}$};    

    \filldraw [black] (1.3,1) circle (1pt) node[anchor=south]{$(u,1)$};
    \filldraw [black] (0,0) circle (1pt) node[anchor=north]{$(0,0)$};
    \filldraw [black] (1,0.76) circle (1pt) node[anchor=north west]{$\Tilde{\mathcal{P}}_{j^*}(u,1)$};
    \filldraw [black] (-0.707,0.707) circle (0.1pt) node[anchor=east]{$\mathbb{S}_+$};
    \filldraw [black] (1,0) circle (1pt) node[anchor=north]{$(\varphi_{j^*},0)$};

    \draw[-] (0,0) -- (1.3,1);
    \end{tikzpicture}
\caption{Projection $\Tilde{\mathcal{P}}_{j^*}$ from phase-space $X^\alpha\hookrightarrow L^2\times \{ 1\}$ into the hyperplanes $C_{j^*}$, which is tangent (at the point $(\varphi_{j^*},0)$) to the sphere at infinity $\mathcal{S}^\infty$.}\label{fig:Ptilde}
\end{figure}

Note the plane $C_{j^*}$ can be rewritten in the coordinates $(\xi,\zeta)$ as
\begin{equation}
    C_{j^*}:=\{ (\xi,\zeta)\in L^2\times \mathbb{R}\text{ $|$ } \xi_{j^*}=+1,\xi_j\in\mathbb{R}\text{ for all } j\in\mathbb{N}_0 \backslash\{j^*\}\}.
\end{equation}
The projection $\tilde{\mathcal{P}}_{j^*}$ induces a semiflow on $C_{j^*}$ through differentiation of \eqref{Ptilde} with respect to time, yielding 
\begin{subequations}\label{flowCj*}
\begin{align}
    \xi_t&
    =\mathcal{L}_\zeta(\xi)-\langle \mathcal{L}_\zeta(\xi),\varphi_{j^*} \rangle \xi  \label{xit}\\
    \zeta_t&
    =-\langle \mathcal{L}_\zeta(\xi), \varphi_{j^*} \rangle \zeta, \label{zetat}
\end{align}
\end{subequations}
where the projected vector field is a homothety of the original vector field \eqref{PDE} with scale factor $\zeta:=\langle u,\varphi_{j^*} \rangle^{-1}$, i.e. $\mathcal{L}_\zeta(\xi):=\zeta \mathcal{L}(\zeta^{-1}\xi)$, as in \eqref{homotvf}.

We are now able to describe the dynamics within the sphere at infinity $\mathcal{S}^\infty$.
\begin{prop}\label{lem:unbddeq}
    There are $N^\infty:=\lfloor \sqrt{b/a^\infty}\rfloor<\infty$ equilibria at infinity, denoted by $\mathcal{E}^{\infty} = \{\pm\Phi_{j}: j=0,\ldots, N^\infty\} \subseteq \mathcal{S}^\infty$, such that $z(\pm\Phi_{j})=j$. Moreover, 
    there is a heteroclinic at infinity, in $\mathcal{H}^\infty \subseteq \mathcal{S}^\infty$, from $\Phi_{j}$ to $\Phi_{k}$ as in \eqref{HetsInfty} if, and only if, $j>k$.
\end{prop}
\begin{proof}
The proof is divided in two parts: first, we construct the equilibria at infinity $\mathcal{E}^\infty\subseteq\mathcal{S}^\infty$; second, we prove the existence of heteroclinics at infinity $\mathcal{H}^\infty\subseteq\mathcal{S}^\infty$.

Let us describe the objects $ \{\pm \Phi_{j} \}_{j=0}^{N^\infty}$ and show they play the role of equilibria at infinity. 
In order to approximate the semiflow at $\mathcal{S}^\infty$, consider $\zeta\rightarrow 0$, which yields the tangent space of $\mathcal{S}^\infty$ at the point $(\varphi_{j^*},0)$, given by $C_{j^*}|_{\zeta=0}$. This tangent space is invariant, since the semiflow in equation \eqref{flowCj*} yields $\zeta_t= 0$.
Indeed, this occurs due to boundedness of the inner product in right hand side of equation \eqref{zetat}:
\begin{align}
    |\langle \mathcal{L}_\zeta(\xi),\varphi_{j^*} \rangle | &\leq
    a|\langle\xi_{xx}+b\xi+f_\zeta(\xi),\varphi_{j^*}\rangle | \nonumber\\
    &\leq a|\lambda_{j^*}|+b+|\langle f_\zeta(\xi),\varphi_{j^*}\rangle|<\infty,
\end{align}
since $a(x,\chi,\chi_x)$ is positive and bounded by some $a\in\mathbb{R}_+$, and $\xi_{j^*}=1$ in $C_{j^*}$. 

In order to dissect the induced semiflow of equation \eqref{xit} in $C_{j^*}$, we write each $\varphi_j$-eigendirectional semiflow as $\xi_j:=\langle \xi ,\varphi_j\rangle$, which satisfies
\begin{align}\label{xij}
   (\xi_j)_t
   = \langle &a_\zeta(x,\xi,\xi_x)\xi_{xx}+b\xi+f_\zeta(x,\xi,\xi_x),\varphi_j \rangle \nonumber\\ - &\langle  a_\zeta(x,\xi,\xi_x)\xi_{xx}+b\xi+f_\zeta(x,\xi,\xi_x),\varphi_{j^*}  \rangle\xi_j,
\end{align}
and then compute the limit as $\zeta\rightarrow 0$, noticing $f_\zeta(\xi)\to 0$, 
$\xi_{j^*}=1$ in $C_{j^*}$, and \eqref{alim}:
\begin{align}\label{xijlim}
   (\xi_j)_t
   &= \langle a^\infty\xi_{xx},\varphi_j \rangle- \langle  a^\infty\xi_{xx},\varphi_{j^*}  \rangle\xi_j\nonumber\\ 
   &=a^\infty(\lambda_{j}-\lambda_{{j^*}})\xi_j.
\end{align}
%
%
%
%
Therefore, grow-up behavior of the solutions $u(t)$ induces the linear flow \eqref{xijlim} in the projected coordinates $(\xi,0)\in C_{j^*}|_{\zeta=0}$. In particular, equilibria of \eqref{xijlim} occur when $\xi_{j}=0$ for all $j\in\mathbb{N}_0$, except $\xi_{j^*}\neq 0$. The only of these equilibria that lie on the sphere at infinity $\mathcal{S}^\infty$ is when $\xi_{j^*}=\pm 1$. Thus, the equilibria at infinity are given by
\begin{equation}\label{EQcsi}
    \pm\Phi_{j^*} := \{ (\xi,0)\in \mathcal{S}^\infty:\xi_{j^*}=\pm 1, \mbox{ and } \xi_{j}=0 \ \forall j \neq j^* \},
\end{equation}
for all $j^*\in \mathbb{N}_{0}$. 

The coordinates $(\chi,z)$ and $(\xi,\zeta)$ are related by means of colinearity,
\begin{equation}\label{changeofcoord}
    (\xi,\zeta)=\frac{1}{\langle \chi , \varphi_{j^*}\rangle}(\chi,z).
\end{equation}
In particular, the equilibria at $\mathcal{S}^\infty$ in terms of $(\xi,\zeta)$-coordinates, given by \eqref{EQcsi}, can be translated to $(\chi,z)$-coordinates:
\begin{equation}\label{EQchi}
\pm\Phi_{j^*} = \{ (\chi,0)\in \mathcal{S}^\infty:\chi_{j^*}=\pm 1, \mbox{ and } \chi_{j}=0 \ \forall j \neq j^* \} \ ,
\end{equation}
%
%
%
%
%
We label the equilibria by $\pm\Phi_{j^*}$ because they \emph{are} the eigenfunction $(\pm\varphi_{j^*},0)\in\mathcal{S}^\infty$:
\begin{equation}
    \pm\Phi_{j^*}=(\pm\varphi_{j^*},0).
\end{equation}
Indeed, recall the equilibria $\pm\Phi_{j^*}$ are described in \eqref{EQchi} by $\chi_{j^*}=\langle \chi, \varphi_{j^*} \rangle=\pm 1$, and all others $\chi_{j}=\langle \chi, \varphi_{j^*} \rangle=0$ for $j\neq j^*$. This means that $\chi=\pm\varphi_{j^*}$, since $\chi$ must lie in the unitary sphere in $L^2$, i.e., $(\chi,0)\in\mathcal{S}^\infty$.
Thus, the zero numbers of the equilibria $\pm\Phi_{k}$, and the difference of equilibria at infinity $\Phi_{k}-\Phi_{j}$, are well defined, since they respectively correspond to the eigenfunctions $\pm\varphi_{k}$ and $\varphi_{k}-\varphi_{j}$.

Let us now construct the heteroclinic network at infinity, $\mathcal{H}^\infty\subseteq \mathcal{S}^\infty$. 
Indeed, given an equilibrium $\Phi_j$ with $j\in \{1,...,N^\infty\}$, we prove that there is a heteroclinic connection from $\pm\Phi_{j}$ to $\pm\Phi_{k}$ for each $k\in\{0,1,\cdots,j-1\}$.

%
%
Consider the hyperplane $C_{j}$, where $\zeta=0$, which contains the equilibrium $\pm\Phi_{j}$.
Then, due to \eqref{xijlim}, each coordinate $\xi_k(t)$ belongs to a linear subspace where its evolution is given by 
\begin{equation}
    (\xi_k)_t=a^\infty(\lambda_{k}-\lambda_{j})\xi_k.
\end{equation}
Thus, $\xi_k(t)$ expands outside $\pm\Phi_{j}$, since $\lambda_{k}>\lambda_{j}$. Moreover, $\xi_k(t_*)=1$ for some $t_*>0$, i.e., $\xi_k(t_*)$ intersects the plane $C_{k}$. Denote such intersection point by $p_*$.

On the other hand, the evolution of $\xi_j(t)$ in the plane $C_{k}$ restricted to $\zeta=0$ yields a contraction in the $\xi_j$-direction of the equilibria $\pm\Phi_{k}\in C_{k}$, with flow given by
\begin{equation}
    (\xi_j)_t=a^\infty(\lambda_{j}-\lambda_{k})\xi_j,
\end{equation}
since $\lambda_{k}>\lambda_{j}$. In particular, contraction occurs for the initial data $p_*$ defined above.

Lastly, note that the expansion and contraction occur in the $C_{j}$ and $C_{k}$ planes, respectively, when $\zeta=0$. Moreover, those are the respective tangent spaces of the sphere at infinity, $\mathcal{S}^\infty$, at the points $\Phi_{j}$ and $\Phi_{k}$. By means of \eqref{changeofcoord}, with $z=0$ and $\zeta=0$, one obtain a topological equivalence of the flow between these spaces. Therefore, we obtain the heteroclinic network at infinity in the sphere at infinity, $\mathcal{S}^\infty$. See Figure \ref{FIG:last}.
\end{proof}
\begin{figure}[H]\centering
    \begin{tikzpicture}[scale=2]
    \draw[->] (-0.1,0) -- (1.5,0) node[right] {$\xi_k$};
    \draw[->] (0,-0.1) -- (0,1.5) node[above] {$\xi_j$};

    \draw[-] (-0.1,1) -- (1.5,1) node[right] {$C_{j}|_{\zeta=0}$};
    \draw[-] (1,-0.1) -- (1,1.5) node[above] {$C_{k}|_{\zeta=0}$};    

    \draw [very thick, domain=0:1.57,variable=\t,smooth] plot ({cos(\t r)},{sin(\t r)}); 
    
    \draw[->,thick,>=stealth] (0.49,1) -- (0.51,1);
    \draw[->,thick,>=stealth] (1,0.51) -- (1,0.49);
    \draw [<-,very thick,>=stealth,domain=0.78:0.79, variable=\t,smooth] plot ({cos(\t r)},{sin(\t r)}) node[anchor=north east] {$\mathcal{S}^\infty$}; 

    \filldraw [black] (1,0) circle (1pt) node[anchor=north west]{$\Phi_{k}$};
    \filldraw [black] (0,1) circle (1pt) node[anchor= south east]{$\Phi_{j}$};
    \filldraw [black] (1,1) circle (1pt) node[anchor= south west]{$p_{*}$};

    \end{tikzpicture}
\caption{The heteroclinics at infinity, $\mathcal{H}^\infty\subseteq \mathcal{S}^\infty$, between two equilibria $\Phi_{j}$ and $\Phi_{k}$ with $j>k$, can be obtained by means of the linear flows in $C_j|_{\zeta=0}$ and $C_k|_{\zeta=0}$.}\label{FIG:last}
\end{figure}

%

\subsection{Grow-up behavior: attainable sphere at infinity}\label{sec:growup}

Next we address the existence of grow-up solutions, $\mathcal{H}^{up}$, which can be seen as heteroclinics from bounded equilibria to equilibria at infinity. First, we prove that the number of zeros of grow-up solutions do not decrease in the limit $t\to\infty$. In particular, grow-up solutions $u(t)$ converge to an equilibrium at infinity with the same number of zeros as $u(t)$, for sufficiently large times. Later, we use the $y$-map to show the blocking and liberalism principles for grow-up solutions, which dictate how heteroclinics can be blocked, and how heteroclinics exist if blocking does not occur, respectively.
These principles yield item 2 in Theorem \ref{attractorthm}.

\begin{lem}\label{NoZeroDrop}
Let $u(t)$ be a grow-up solution in the unstable manifold of a bounded equilibrium $e_j\in \mathcal{E}^b$ such that
\begin{equation}\label{saopaulo2018}
    z(u(t)-e_j)=k, \qquad \iota :=sign(u(t,0)-e_j(0))\in \{\pm \},
\end{equation}
for sufficiently large time $t$. Then, the projected solution $\mathcal{P}(u(t))$ converges (in $C^1$) to the unbounded equilibrium $\iota\Phi_{k}\in\mathcal{E}^\infty$.
\end{lem} 
\begin{proof}
Due to \eqref{growupdirection}, the normalization of the grow-up solution $u(t)$ either converges to $+\Phi_{l}$, or to $-\Phi_l$, namely
\begin{equation}\label{L2conv}
    \lim_{t\rightarrow \infty} \left\Vert  \frac{u(t)}{\Vert u(t) \Vert_{L^2}} -\varphi_{l} \right\Vert_{L^2} =0, \quad \text{ or } \quad \lim_{t\rightarrow \infty} \left\Vert  \frac{u(t)}{\Vert u(t) \Vert_{L^2}} +\varphi_{l} \right\Vert_{L^2} =0.
\end{equation}
Without loss of generality, we consider the case $\iota=+$, and thus comparison implies that $u(t,0)>e_j(0)$ for all sufficiently large $t>0$. This implies that $\mathcal{P}(u(t))$ converges in $L^2$ to the equilibrium $+\Phi_{l}$ with positive sign.

It is enough to show that the convergence \eqref{L2conv} also takes place in the $C^1$ topology.
Indeed, if the convergence occurs in $C^1$, then the limit of $u(t)$ has a constant number of zeros for sufficiently large time $t$, given by $z(u(t))=z(u(t)-e_j)$, and does not drop at $t=\infty$. Thus, hypothesis \eqref{saopaulo2018} implies $l=k$.
    
We rewrite $u(t)=p(t)+q(t)$, as in \eqref{ultimodia}, where $q(t)$ is bounded. As a consequence, \eqref{L2conv} implies that
\begin{equation}\label{pL2conv}
    \lim_{t\rightarrow \infty} \left\Vert  \frac{p(t)}{\Vert u(t) \Vert_{L^2}} -\varphi_{l} \right\Vert_{L^2} =0,
\end{equation}
since $u(t)$ is a grow-up solution, and hence $q(t)/\Vert u(t) \Vert_{L^2}\to 0$.

Recall $p(t)$ is contained in a finite dimensional subspace, as well as the span of $\varphi_{l}^\pm $ which is one dimensional, and thus norms in such finite dimensional subspace are equivalent. In particular, convergence in \eqref{pL2conv} also holds in terms of the $C^1$-topology. 
Therefore, the fact that $q(t)/\Vert u(t) \Vert_{L^2}\to 0$ implies that the convergence \eqref{L2conv} happens in $C^1$.
\end{proof}
%
%
%
\begin{prop}\emph{\textbf{Blocking at infinity.} }\label{prop:IB}
If the bounded equilibrium, $e_j\in \mathcal{E}$, and the equilibrium at infinity, $\pm\Phi_{k}\in \mathcal{E}^\infty$, are not adjacent, then there does not exist a grow-up solution, which is a heteroclinic in $\mathcal{H}^{up}$, that connects them according to \eqref{Hetup}.
\end{prop}
\begin{proof}
Assume, towards a contradiction, the existence of a grow-up solution $u(t)$ which is contained in the unstable manifold of $e_j$ and converges (in forward time) to $\Phi_{k}$. We know that any grow-up solution $u(t)\in W^u(e_j)$ converges to an equilibria at infinity, $\{ \pm\Phi_{l}:l,\ldots,N^\infty \}$. Since $e_j$ and $\Phi_{k}$ are not adjacent, there exists a bounded equilibrium $e_*\in \mathcal{E}$ satisfying
$e_j(0)<e_*(0)<\Phi_{k}(0)$ and $z(e_j-e_*)=z(\Phi_{k}-e_*)=z(\Phi_{k}-e_j)$.

Due to the non-adjacency and the $C^1$-convergence in Lemma \ref{NoZeroDrop}, that
\begin{equation}\label{1}
    z(u(-t)-e_*)=z(e_j-e_*)=z(\Phi_{k}-e_*)=z(u(t)-e_*)
\end{equation}
for sufficiently large $t$. 

On the other hand, note that $e_j(0)<e_*(0)<u(t,0)$ for sufficiently large values of $t$, and therefore $e_j(0)-e_*(0)<0<u(t,0)-e_*(0)$, i.e., $e_j(0)-e_*(0)$ and $u(t,0)-e_*(0)$ have opposite signs. Then, the solution $\Tilde{u}(t):=u(t)-e_*$ has a multiple zero in the boundary, due to Neumann boundary conditions. Consequently, there exists a large dropping time in the boundary, due to the Lemma \ref{droplem}. Therefore, 
\begin{equation}\label{2}
    z(u(-t)-e_*)>z(u(t)-e_*)    
\end{equation}
for sufficiently large $t$.
Equations \eqref{1} and \eqref{2} yield a contradiction.
\end{proof}
\begin{prop}\emph{\textbf{Liberalism at infinity.} }\label{prop:IL}
If the bounded equilibrium, $e_j\in \mathcal{E}$, and the equilibrium at infinity, $\pm\Phi_{k}\in \mathcal{E}^\infty$, are adjacent, then there exists a grow-up solution, which is a heteroclinic in $\mathcal{H}^{up}$, that connects them according to \eqref{Hetup}.
\end{prop}
\begin{proof}
Corollary \ref{ymapcor} guarantees the existence of a solution $u(t)\in W^u(e_j)$ satisfying
\begin{equation}\label{constz}
z(u(t)- e_j)=k,
\end{equation}
for all $t\geq0$, so that $\iota:=\text{sign}(u(t,0)-e_j(0))\in\{\pm\}$ is indeed constant for all $t\geq 0$.

Note that $u(t)$ grows up. Indeed, suppose towards a contradiction that it converges to a bounded equilibrium, i.e. $\lim_{t\rightarrow \infty}u(t)=e_*\in \mathcal{E}$. The Lemma \ref{droplem} and equation \eqref{constz} imply that $z(e_*-e_j)\leq k$. But adjacency prevents the equality, since in addition to $z(e_*-e_j)= k$, we would have that $\text{sign} (e_*(0)-e_j(0))=\pm$. Therefore, such $e_*$ necessarily satisfies $z(e_*-e_j)< k$, which implies that the zero number of the shifted solution $u(t,0)-e_j(0)$ has to drop at $t=\infty$. This cannot happen, since $e_*-e_j$ only has simple zeros, and hence $z(u(t)-e_j)$ would have to drop at some finite time. This contradicts \eqref{constz}, and thereby $u(t)$ grows-up.

Then, convergence of $u(t)$ to $\pm\Phi_{k}$ as $t\rightarrow \infty$ is achieved by Lemma \ref{NoZeroDrop}, which implies that unbounded solutions cannot have zero dropping at $t=\infty$.
\end{proof}

\section{Discussion}\label{sec:disc}

We give three comments regarding our current results. First on the hypothesis \eqref{alim}, second on the possibility of using our tools to obtain results for geometric PDEs, and third on the applicability to construct initial data for certain black holes.

Not all quasilinear diffusion coefficients admit a well-defined limiting semiflow at $\mathcal{S}^\infty$, since they may oscillate or diverge at infinity. We have assumed that the diffusion coefficient is asymptotically constant in \eqref{alim} to tame such a problem.
We believe that our methods can be extended to a broader class of diffusion with a well-defined limit in the sphere at infinity, e.g. $\lim_{z\to 0} a^z\left(x,\chi,\chi_x\right)=a^\infty(x,\chi,\chi_x)\geq \epsilon >0$.
However, such limiting diffusion coefficient yields a highly nonlinear and nonlocal induced semiflow at $\mathcal{S}^\infty$.
Indeed, the induced PDE at infinity in \eqref{flowSPHEREz0} would not be semilinear anymore, but quasilinear, and the projected equation \eqref{xijlim} would not be linear. 
These type of limiting diffusion may support a different behavior within the sphere at infinity, in contrast to the Chafee-Infante structure we have constructed.
In particular, we conjecture that any (bounded) Sturm attractor may be realized  as an extended unbounded attractor in \eqref{extendedA} for some $a$ and $f$, containing a prescribed dynamics at the sphere at infinity, see \cite{FiedlerRocha99,FiedlerRocha00}.



Next, note that the inverse mean curvature flow (IMCF) for hypersurfaces yields a fully nonlinear parabolic equation with grow-up behavior for a broad class of initial data.
In particular, the evolution of star shaped (or more general) initial data under the rescaled ICMF converges to a round sphere, yielding certain stability result for the round sphere, see \cite{Gerhardt,Harvie}. 
In our approach, such a stability result can be interpreted as follows: the compactification of the star shaped (or more general grow-up) subset of phase-space of IMCF should have a unique stable equilibrium at infinity, which corresponds to the round sphere. 
Our present results can be extended for fully nonlinear equations with an induced gradient semiflow at infinity, akin to \cite{LappicyFully,LappicyFiedler18,LappicyBeatriz}, under reasonable hypothesis on $f$. Hence, our compactification procedure provides a method that may prove stability of a broader subset of initial data for certain geometric flows that admit grow-up. 


Lastly, consider a variant of \eqref{PDE} with a singular diffusion coefficient, given by
\begin{equation}\label{PDEsing}
    u_t=a(x,u,u_x)\left[u_{xx}+\frac{u_x}{\tan(x)}\right]+f(x,u,u_x),
\end{equation}
where $x\in [0,\pi]$ has Neumann boundary conditions.
This equation describes the space of initial data for certain axially symmetric black holes, where $t$ is a rescaled radial distance from the singularity with event horizon at $t\to\infty$, $u$ is a component of a Riemannian metric for a spacelike hypersurface, and $f$ is related to a prescribed scalar curvature of said Riemannian manifold. See \cite{Bartnik,Smith09,FiedlerHellSmith,LappicyBlackHoles}.

In contrast to \cite{Bartnik,Smith09}, who have freely specified the metric coefficient at the event horizon to construct exterior metrics, the authors in \cite{FiedlerHellSmith} have shown that black holes with same exterior metric admit a plethora of different metrics inside the horizon, by means of equivariant bifurcation.
Moreover, some metric coefficients can not be freely prescribed at the event horizon, but it should be constrained to the equilibria within the unbounded attractor of equation \eqref{PDEsing}, see \cite{LappicyBlackHoles}. 
We mention that the bounded equilibria and equilibria at infinity have different interpretations: the former correspond to self-similar interior Schwarzschild solutions, whereas the latter amount to non-self-similar interior Schwarzschild solutions, 
see \cite{LappicyBlackHoles}.
Our current results can be extended to the case of a singular diffusion, as in \eqref{PDEsing}, according to \cite{LappicySing}. However, in order to fully understand the structure at the event horizon, one has to construct the unbounded attractor for the degenerate diffusion $a(u)=u^2$.

\textbf{Acknowledgment.} PL was firstly supported by FAPESP, 17/07882-0, and later on by Marie Skłodowska-Curie Actions, UNA4CAREER H2020 Cofund, 847635, 
with the project DYNCOSMOS. JF was supported by FAPESP, 16/04925-7, and partially supported by CAPES - Programa CAPES-PrInt, Processo No. 88881.311616/2018-00. 


\medskip

\begin{thebibliography}{10}
\small{

\bibitem{Amann95}
H.~Amann.
Linear and Quasilinear Parabolic Problems, Volume I
\emph{Birkh{\"a}user Basel} (1988).

\bibitem{Angenent88}
S.~Angenent.
The zero set of a solution of a parabolic equation.
\emph{J. f{\"u}r die reine und angewandte Math.} \textbf{390}, 79 -- 96, (1988).

\bibitem{ArrietaCarvalhoBernal00}
J.~Arrieta, A.N.~Carvalho, A.~Rodríguez-Bernal.
Attractors of parabolic problems with nonlinear boundary conditions. uniform bounds
\emph{Comm. P.D.E.} \textbf{25}, 1 -- 37, (2000).

\bibitem{ArrietaPardo}
J.~Arrieta, R.~Pardo and A. Rodríguez-Bernal.
Asymptotic behavior of degenerate logistic equations.
\emph{J. Diff. Eq.} \textbf{259}, 6368 -- 6398, (2015).

\bibitem{BabinVishik92}
A.V.~Babin and M.I.~Vishik.
Attractors of Evolution Equations.
\emph{Elsevier Science}, (1992).

\bibitem{Bartnik}
R.~Bartnik.
Quasi-spherical metrics and prescribed scalar curvature.
\emph{J. Diff. Geometry} \textbf{37}, 31 -- 71, (1993).

\bibitem{BenGal10}
N.~Ben-Gal.
Grow-Up Solutions and Heteroclinics to Infinity for Scalar Parabolic PDEs.
\emph{Ph.D. Thesis, Division of Applied Mathematics, Brown University}, (2010).

\bibitem{BortolanJu22}
M.~Bortolan and J.~Fernandes.
Sufficient Conditions for the Existence and Uniqueness of Maximal Attractors for Autonomous and Nonautonomous Dynamical Systems.
\emph{J. Dyn. Diff. Eq.} \textbf{1}, 1 -- 30, (2022).

\bibitem{BrunovskyFiedler86}
P.~Brunovsk{\'y} and B.~Fiedler.
Numbers of Zeros on Invariant Manifolds in Reaction-diffusion Equations
\emph{Nonlinear Analysis: TMA} \textbf{10}, 179--193, (1986).

\bibitem{FiedlerBrunovsky88}
P.~Brunovsk{\'y} and B.~Fiedler.
Connecting orbits in scalar reaction diffusion equations.
\emph{Dynamics Reported} \textbf{1}, 57--89, (1988).

\bibitem{FiedlerBrunovsky89}
P.~Brunovsk{\'y} and B.~Fiedler.
Connecting orbits in scalar reaction diffusion equations II: The
  complete solution.
\emph{J. Diff. Eq.} \textbf{81}, 106--135, (1989).

\bibitem{BCP}
S. Bruschi, A. N. Carvalho, J. Pimentel.
Limiting grow-up behavior for a one-parameter family of dissipative PDEs. 
\emph{Indiana Univ. Math. J.} \textbf{69}, 657--683, (2020)

\bibitem{CG92}
V.~Chepyzhov and A.~Goritskii. 
Unbounded attractors of evolution equations.
\emph{Properties of Global Attractors of Partial Differential Equations. Adv. Soviet
Math. 10, Amer. Math. Soc., Providence, RI}, eds. A. V. Babin and M. I. Vishik, 85--128, (1992)

\bibitem{FerreiradePablo20}
R.~Ferreira and A.~de Pablo.
Grow-up for a quasilinear heat equation with a localized reaction.
 \emph{J. Diff. Eq.} \textbf{268}, 6211--6229, (2020).

\bibitem{FiedlerHellSmith}
B.~Fiedler, J.~Hell, and B.~Smith.
Anisotropic Einstein data with isotropic nonnegative scalar curvature.
\emph{Ann. Inst. Henri Poincar{\'e} - An. Non Lin\'eaire} \textbf{32}, 401--428, (2015).

\bibitem{FiedlerRocha96}
B.~Fiedler and C.~Rocha.
Heteroclinic orbits of semilinear parabolic equations.
\emph{J. Diff. Eq.} \textbf{125}, 239--281, (1996).

\bibitem{FiedlerRocha99}
B.~Fiedler, C.~Rocha.
Realization of meander permutations by boundary value problems.
\emph{J. Diff. Eq.} \textbf{156}, 282 -- 308, (1999).

\bibitem{FiedlerRocha00}
B.~Fiedler and C.~Rocha.
Orbit Equivalence of Global Attractors of Semilinear Parabolic Differential Equations.
\emph{Trans. Am. Math. Soc.} \textbf{352}, 257--284, (2000).

\bibitem{FRW}
B.~Fiedler, C.~Rocha and M.~Wolfrum.
Heteroclinic orbits of semilinear parabolic equations.
\emph{J. Diff. Eq.} \textbf{201}, 99--138, (2004).

\bibitem{FiMa02}
M.~Fila and H.~Matano.
Blow-up in Nonlinear Heat Equations from the Dynamical Systems Point of View.
\emph{Handbook of Dyn. sys., vol. 2,} ed. B. Fiedler,723--758, (2002).

\bibitem{FuscoRocha}
G.~Fusco and C.~Rocha.
A permutation related to the dynamics of a scalar parabolic PDE.
\emph{J. Diff. Eq.} \textbf{91}, 111--137, (1991).

\bibitem{Harvie}
B.~Harvie.
Inverse Mean Curvature Flow over Non-Star-Shaped Surfaces.
\emph{Math. Res. Letters} \textbf{29}, 1065--1086, (2022). 

\bibitem{Hell09}
J.~Hell.
Conley Index at Infinity. 
\emph{Topol. Methods Nonlin. Anal.} \textbf{42}, 137--167, (2013).

\bibitem{Gerhardt}
C.~Gerhardt.
Flow of nonconvex hypersurfaces into spheres.
\emph{J. Diff. Geom.} \textbf{32}, 299--314, (1990).

\bibitem{Henry81}
D.~Henry.
\emph {Geometric Theory of Semilinear Parabolic Equations}.
Springer-Verlag New York, (1981).

\bibitem{LappicyQuasi}
P.~Lappicy.
Sturm attractors for quasilinear parabolic equations.
\emph{J. Diff. Eq.} \textbf{265}, 4642--4660, (2018).

\bibitem{LappicyBlackHoles}
P.~Lappicy.
Space of initial data for self-similar Schwarzschild solutions of the Einstein equations. 
\emph{Phys. Rev. D} \textbf{99}, 043509, (2019).

\bibitem{LappicySing}
P.~Lappicy.
Sturm attractors for quasilinear parabolic equations with singular coefficients.
\emph{J. Dyn. Diff. Eq.}  \textbf{32}, 359--390, (2020).

\bibitem{LappicyFiedler18}
P.~Lappicy and B.~Fiedler.
A Lyapunov function for fully nonlinear parabolic equations in one spatial variable.
\emph{São Paulo J. Math. Sc.} \textbf{13}, 283--291, (2019).

\bibitem{LappicyFully}
P.~Lappicy.
Sturm attractors for fully nonlinear parabolic equations in one spatial dimension.
\emph{Rev. Mat. Complutense} \textbf{36}, 725--747, (2023).

\bibitem{LappicyBeatriz}
P.~Lappicy and E.~Beatriz.
An energy formula for fully nonlinear degenerate parabolic equations in one spatial dimension.
\emph{To appear in Math. Annalen, https://doi.org/10.1007/s00208-023-02740-5}, (2023).

\bibitem{LopezGomez01}
J.~L\'opez-G\'omez.
Approaching metasolutions by solutions.
\emph{Diff. Integral Eq.} \textbf{14}, 739--750, (2001).

\bibitem{LopezGomez03}
J.~L\'opez-G\'omez.
Dynamics of parabolic equations: from classical solutions to metasolutions.
\emph{Diff. Integral Eq.} \textbf{16}, 813--828, (2003).

\bibitem{LopezGomez16}
J.~L\'opez-G\'omez.
\emph{Metasolutions of parabolic equations in population dynamics.}
Boca Raton, FL: CRC Press, (2016).

\bibitem{Lunardi95}
A.~Lunardi.
\emph{Analytic Semigroups and Optimal Regularity in Parabolic Problems}.
Springer Basel, (1995).

\bibitem{Matano82}
H.~Matano.
Non increase of the lapnumber for a one dimensional semilinear parabolic equation.
\emph{J. Fac. Sci. Univ. Tokyo IA Math} \textbf{29}, 401--441, (1982).
  
\bibitem{Matano88}
H.~Matano.
Asymptotic behavior of solutions of semilinear heat equations on
  $S^1$.
\emph{Nonlinear Diffusion Equations and Their Equilibrium States II,}
  eds. W.-M. Ni, L. A. Peletier, J. Serrin, 139--162, (1988).

\bibitem{Mik}
M.~Miklavčič.
A sharp condition for existence of an inertial manifold.
\emph{J. Dyn. Diff. Eq.} \textbf{3}, 437--456, (1991).

\bibitem{MussoSire19}
M.~Musso, Y.~Sire, J.~Wei, Y.~Zheng and Y.~Zhou .
Infinite time blow-up for the fractional heat equation with critical exponent.
\emph{Math. Annalen} \textbf{375}, 361--424, (2019).

\bibitem{P}
J.~Pimentel.
Unbounded Sturm global attractors for semilinear parabolic equations on the circle. 
\emph{SIAM Journal on Mathematical Analysis} \textbf{48} 3860 -- 3882, (2016).

\bibitem{RochaPimentel16}
J.~Pimentel and C.~Rocha.
A permutation related to non-compact global attractors for slowly non-dissipative systems.
\emph{J. Dyn. Diff. Eq.} \textbf{28}, 1--15, (2016).


\bibitem{QuittnerSouplet19}
P.~Quittner and P.~Souplet.
\emph{Superlinear Parabolic Problems: Blow-up, Global Existence and Steady States.} Birkhäuser, 2nd Ed., (2019).



\bibitem{Smith09}
B.~Smith.
Black hole initial data with a horizon of prescribed geometry.
\emph{Gen. Relativ. Gravit.} \textbf{41}, 1013 -- 1024, (2009).

\bibitem{Sturm}
C.~Sturm.
Sur une classe d'{\'e}quations {\`a} diff{\'e}rences partielles.
 \emph{J. Math. Pures. Appl. I} \textbf{1}, 373--444, (1836).
 
\bibitem{Wolfrum02}
M.~Wolfrum.
A Sequence of Order Relations: Encoding Heteroclinic Connections in Scalar Parabolic PDE.
 \emph{J. Diff. Eq.} \textbf{183}, 56--78, (2002).

\bibitem{Zelenyak68}
T.~I. Zelenyak.
Stabilization of solutions of boundary value problems for a second
  order parabolic equation with one space variable.
 \emph{Differ. Uravn.} \textbf{4}, 34--45, (1968).
 
}
\end{thebibliography}
\end{document}